\RequirePackage{fix-cm}
\documentclass[smallextended]{svjour3}     
\usepackage{graphicx}
\usepackage{amsmath}
\usepackage{bm}
\usepackage{tabularx}
\usepackage[normalem]{ulem}
\usepackage{ragged2e}
\usepackage{amssymb,amsfonts}
\usepackage{color}
\usepackage{xcolor}
\usepackage[sort]{cite}
\usepackage{enumerate}   

\newcolumntype{L}[1]{>{\raggedright\arraybackslash}p{#1}}
\newcolumntype{C}[1]{>{\centering\arraybackslash}p{#1}}
\newcolumntype{R}[1]{>{\raggedleft\arraybackslash}p{#1}}
\newcolumntype{J}[1]{>{\justifying\arraybackslash}p{#1}}

\newcommand{\code}[1]{\texttt{#1}}
\newcommand{\ten}[1]{\ensuremath{\boldsymbol{\mathsf{#1}}}}
\newcommand{\tI}{{\ensuremath{\ten I}}}
\renewcommand{\vec}[1]{{\ensuremath{\boldsymbol{\mathrm #1}}}}

\newcommand{\FF}{\mathrm{ff}}
\newcommand{\PM}{\mathrm{pm}}

\newcommand{\Pt}{\mathcal{P}_{\mathrm{{tri}}}}
\newcommand{\Pc}{\mathcal{P}_{\mathrm{{con}}}}
\newcommand{\Pd}{\mathcal{P}_{\mathrm{{diag}}}}

\newcommand{\hPd}{\hat{\mathcal{P}}_{\mathrm{{diag}}}}
\newcommand{\hPt}{\hat{\mathcal{P}}_{\mathrm{{tri}}}}
\newcommand{\hPc}{\hat{\mathcal{P}}_{\mathrm{{con}}}}
\newcommand{\Abar}{\overline{\mathcal{A}}}

\newcommand{\ABJS}{\mathcal{A}_{\mathrm{BJS}}}
\newcommand{\ABJ}{\mathcal{A}_{\mathrm{BJ}}}

\newcommand{\ex}{\mathrm{ex}}

\begin{document}

\title{Efficient preconditioners for coupled Stokes--Darcy problems with MAC scheme\thanks{
The work is funded by the Deutsche Forschungsgemeinschaft (DFG, German Research Foundation) – Project Number 327154368 – SFB 1313 and Project Number 490872182.
}
}
\subtitle{Spectral analysis and numerical study}

\titlerunning{Preconditioners for coupled Stokes--Darcy problems}        

\author{Paula Strohbeck \and
        Iryna Rybak
}


\institute{P. Strohbeck \at
              University of Stuttgart, Institute of Applied Analysis and Numerical Simulation, Pfaffenwaldring 57, 70569 Stuttgart, Germany \\
              \email{paula.strohbeck@ians.uni-stuttgart.de}           
           \and
           I. Rybak \at
        University of Stuttgart, Institute of Applied Analysis and Numerical Simulation, Pfaffenwaldring 57, 70569 Stuttgart, Germany \\
        \email{iryna.rybak@ians.uni-stuttgart.de}
}

\date{Received: date / Accepted: date}

\maketitle

\begin{abstract}
Coupled systems of free flow and porous media arise in a variety of technical and environmental applications. For laminar flow regimes, such systems are described by the Stokes equations in the free-flow region and Darcy's law in the porous medium. An appropriate set of coupling conditions is needed on the fluid--porous interface. Discretisations of the Stokes--Darcy problems yield large, sparse, ill-conditioned, and, depending on the interface conditions, non-symmetric linear systems. Therefore, robust and efficient preconditioners are needed to accelerate convergence of the applied Krylov method. In this work, we consider the second order MAC scheme for the coupled Stokes--Darcy problems and develop and investigate block diagonal, block triangular and constraint preconditioners. We apply two classical sets of coupling conditions considering the Beavers--Joseph and the Beavers--Joseph--Saffman condition for the tangential velocity. For the Beavers--Joseph interface condition, the resulting system is non-symmetric, therefore GMRES method is used for both cases. Spectral analysis is conducted for the exact versions of the preconditioners identifying clusters and bounds.
Furthermore, for practical use we develop efficient inexact versions of the preconditioners. We demonstrate effectiveness and robustness of the proposed preconditioners in numerical experiments. 
\end{abstract}

\section{Introduction}
Coupled systems of free flow and porous media appear in a variety of technical applications and environmental settings such as industrial filtration, water/gas management in fuel cells, and surface--subsurface interactions.
For low Reynolds numbers, the free flow is governed by the Stokes equations, and flow in the porous-medium region is described by Darcy's law. Interface conditions are needed to couple these models at the fluid--porous interface. In this work, we consider the classical set of coupling conditions, which consists of the conservation of mass across the interface, the balance of normal forces and either the Beavers--Joseph or the Beavers--Joseph--Saffman interface condition on the tangential velocity, e.g.~\cite{Beavers_Joseph_67,Discacciati-Miglio-Quarteroni-02,Layton-Schieweck-Yotov-03,Saffman}.

Different discretisations for the coupled Stokes--Darcy problems have been investigated such as the finite element method~\cite{Beik_Benzi_22,Cai_etal_09,Chidyagwai_etal_16,Discacciati-Miglio-Quarteroni-02,Layton-Schieweck-Yotov-03}, the finite volume method~\cite{Rybak_etal_15,schmalfuss2021partitioned}, the discontinuous Galerkin method~\cite{Chidyagwai-Riviere-09,Zhou-et-al-19,Zhao-Park-20},  or their combinations. These discretisations yield large, sparse, ill-conditioned and, e.g. in case of the Beavers--Joseph coupling condition, non-symmetric linear systems. The Krylov methods are typically applied to efficiently solve large linear systems. Since the Beavers--Joseph coupling condition leads to non-symmetric matrices and preconditioned systems are non-symmetric, we apply in this paper the GMRES method. Convergence of iterative methods can be significantly enhanced by an appropriate choice of preconditioners~\cite{Benzi_02,Rozloznik,Saad_03,vanderVorst-23}.

Preconditioners for the GMRES method applied to solve the Stokes--Darcy problem discretised by the finite element method have been recently studied~\cite{Beik_Benzi_22,Cai_etal_09,Chidyagwai_etal_16,Boon-et-al-22}. In particular, a block diagonal and a block triangular preconditioner based on decoupling the Stokes--Darcy system were developed in~\cite{Cai_etal_09} for the case of the Beavers--Joseph--Saffman condition on the tangential velocity. A constraint preconditioner for this coupled problem was proposed in~\cite{Chidyagwai_etal_16}. Spectral and field-of-values analysis for the block triangular preconditioner from \cite{Cai_etal_09} and the constraint preconditioners from~\cite{Chidyagwai_etal_16} was conducted in~\cite{Beik_Benzi_22}. 

In this paper, we focus on preconditioners for the Stokes--Darcy system discretised by the finite volume method on staggered grids (MAC scheme), since this discretisation scheme is mass conservative, stable and allows natural coupling across the fluid--porous interface, e.g.~\cite{Eymard_etal_10,Rybak_etal_15,Greif_He_2023}. The discretisation of coupling conditions at the interface is crucial, and in contrast to~\cite{Greif_He_2023,Shiu-Ong-Lai-18} we develop a second order MAC scheme for the coupled problem. The discretisation yields a symmetric matrix $\mathcal{A}_{\mathrm{BJS}}$ for the Beavers--Joseph--Saffman condition and a non-symmetric matrix $\mathcal{A}_{\mathrm{BJ}}$  for the Beavers--Joseph condition. The matrices are displayed in equation~\eqref{eq:A-BJS-BJ} and illustrated in Fig.~\ref{fig:domain-sparsity}. Note that $\mathcal{A}_{\mathrm{BJS}}$ is a double saddle point matrix. There exist several efficient preconditioners for double saddle point problems in the literature, e.g. a block diagonal preconditioner was developed in~\cite{He_etal_21} and several block triangular preconditioners were proposed in~\cite{Beik_Benzi_18,He_etal_21,Huang-et-al-22,Greif_He_2023}. The matrix $\mathcal{A}_{\mathrm{BJS}}$ can also be interpreted as a standard saddle point matrix. This implies that besides results on the double saddle point problems, results on the standard saddle point problems can be applied as well. Different types of preconditioners for such systems were established, e.g.~\cite{Axelsson_Neytcheva_03,Bramble_Pasciak_88,Perugia_Simoncini_00,Rozloznik,Saad_03}.

In this work, we consider three main classes of preconditioners and construct one for each class, namely, one block diagonal, one block triangular, and one constraint preconditioner. In contrast to the work available in the literature, we develop preconditioners not only for the Beavers--Joseph--Saffman interface condition, but also for the more general Beavers--Joseph coupling condition on the tangential velocity.

In this paper, we conduct spectral analysis for the preconditioned systems extending our previous work~\cite{strohbeck2023robust}, where we provided only numerical simulation results and considered other coupling conditions for the Stokes--Darcy problem. 
For the convergence of the GMRES method, evaluating the spectrum is not always sufficient~\cite{Greenbaum_etal_96}. Therefore, the field-of-values analysis is needed~\cite{Benzi_21,Loghin_Wathen_04}. It yields estimates on the residual and guarantees the GMRES convergence independent of the grid width. Such analysis is beyond the scope of this work and forms the basis of a subsequent manuscript.

The direct use of the exact preconditioners is computationally expensive. Thus, accurate and easily invertible approximations are required, e.g.~\cite{Benzi_etal_05,Rozloznik,Saad_03}. In this work, we therefore present also inexact variants of the constructed preconditioners.
Moreover, robust preconditioners are desired such that the choice of the physical parameters does not significantly influence the convergence~\cite{Rodrigo-et-al-23,Antonietti-et-al-20,Boon-et-al-22,Budia_Hu_21,Budia_etal_2022}. We illustrate the robustness and efficiency of the developed preconditioners in numerical experiments.

The paper is structured as follows. In section~\ref{sec:model-formulation}, we present the coupled Stokes--Darcy problems with two sets of interface conditions and introduce the corresponding discrete systems. In section~\ref{sec:preconditioners}, we develop three preconditioners and propose their efficient inexact variants. We provide spectral analysis of the preconditioned systems in section~\ref{sec:spectral-analysis}. To demonstrate efficiency and robustness of the constructed preconditioners, numerical experiments are conducted in section~\ref{sec:numerical-results}. In section~\ref{sec:conclusions}, we summarise the obtained results and present possible extensions of this work. 

\enlargethispage{0.85cm}

\section{Model formulation}
\label{sec:model-formulation}
In this paper, we consider a two-dimensional setting. The coupled flow domain $\Omega = \Omega_{\mathrm{pm}} \cup \Omega_{\mathrm{ff}}\subset \mathbb{R}^2$ consists of the free flow $\Omega_{\mathrm{ff}}$ and the porous medium $\Omega_{\mathrm{pm}}$. The two flow regions are separated by the sharp fluid--porous interface $\Sigma$ (Fig~\ref{fig:domain-sparsity}, left). We consider steady-state, incompressible, single-fluid-phase flows at low Reynolds numbers ($Re \ll 1$). The solid phase is non-deformable and rigid that leads to a constant porosity. We deal with homogeneous isotropic and orthotropic porous media. The whole flow system is assumed to be isothermal.

\subsection{Mathematical model}
Coupled flow formulation consists of two different flow models in the two domains and an appropriate set of coupling conditions on the fluid--porous interface. Under the assumptions on the flow made above, the Stokes equations~\eqref{eq:Stokes} are used in $\Omega_{\mathrm{ff}}$. Without loss of generality, we consider Dirichlet boundary conditions on the external boundary
\begin{align}
\label{eq:Stokes}
        \nabla \cdot \textbf{v}_{\mathrm{ff}} = 0, \qquad - \nabla \cdot \ten{T}(\vec{v}_{\mathrm{ff}},p_{\mathrm{ff}})  &= \vec{f}_{\mathrm{ff}} &&\text{in} \; \; \Omega_{\mathrm{ff}}, \\
        \label{eq:BC-ff}
        \vec{v}_{\mathrm{ff}} &= \overline{\vec{v}} &&\text{on}\; \; \partial \Omega_{\mathrm{ff}} \backslash \Sigma,
\end{align}
where $\textbf{v}_{\mathrm{ff}}$ is the fluid velocity, $p_{\mathrm{ff}}$ is the fluid pressure, $\vec{f}_{\mathrm{ff}}$ is a source term, e.g. body force, and $\overline{\textbf{v}}$ is a given function. We consider the stress tensor $\ten{T}(\vec{v}_{\mathrm{ff}},p_{\mathrm{ff}}) = \mu (\nabla \vec{v}_{\mathrm{ff}} + (\nabla \vec{v}_{\mathrm{ff}})^\top ) - p_{\mathrm{ff}} \tI$, where $\mu$ is the dynamic viscosity and $\tI$ is the identity tensor. 

Fluid flow in the porous-medium domain $\Omega_{\mathrm{pm}}$ is based on Darcy's law. Again, we consider Dirichlet boundary conditions on the external boundary 
\begin{align}
    \label{eq:Darcy-law}
    \nabla \cdot \vec{v}_{\mathrm{pm}} = f_{\mathrm{pm}}, \qquad \vec{v}_{\mathrm{pm}} &= - \mu^{-1}\ten{K}\nabla p_{\mathrm{pm}}&&\text{in} \; \; \Omega_{\mathrm{pm}},\\
    \label{eq:BC-pm}
    p_{\mathrm{pm}} &= \overline{p}&&\text{on}\;\; \partial \Omega_{\mathrm{pm}} \backslash \Sigma,
\end{align}
where $\vec{v}_{\mathrm{pm}}$ is the Darcy velocity, $p_{\mathrm{pm}}$ is the pressure, $f_\PM$ is a source term, $\ten{K}$ is the intrinsic permeability tensor, and $\overline{p}$ is a given function. The permeability tensor is symmetric, positive definite, and bounded. In this paper, we restrict ourselves to isotropic ($\ten{K} = k\,\tI$, $k>0$) and orthotropic ($\ten{K} = \text{diag}\left(k_{xx},k_{yy}\right), \, k_{xx},\, k_{yy} > 0$) porous media.
\vspace{-0.5cm}
\begin{figure}[!ht]
    \centering
    \includegraphics[scale = 0.9]{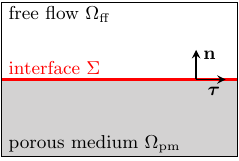} \hspace{5ex} 
     \includegraphics[scale=0.26]{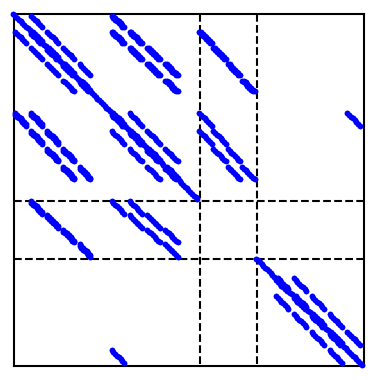}\hspace{5ex} 
     \includegraphics[scale = 0.26]{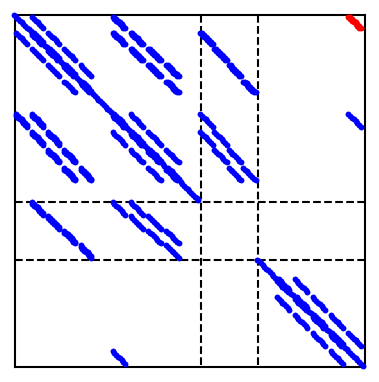}
   \caption{Flow system description (left) and sparsity structures of $\mathcal{A}_{\mathrm{BJS}}$ (middle) and  $\mathcal{A}_{\mathrm{BJ}}$ (right) for $h=1/8$ with symmetric non-zero entries in blue (\textcolor{blue}{$\bullet$}) and non-symmetric non-zero entries in red (\textcolor{red}{$\bullet$})}
    \label{fig:domain-sparsity}
\end{figure}

In addition to the boundary conditions on the external boundary, appropriate coupling conditions have to be defined on the fluid--porous interface $\Sigma$. The classical set of coupling conditions consists of the conservation of mass across the interface~\eqref{eq:IC-mass}, the balance of normal forces~\eqref{eq:IC-momentum} and the Beavers--Joseph condition~\eqref{eq:IC-BJ} on the tangential velocity~\cite{Beavers_Joseph_67}:
\begin{align}
    \vec{v}_{\mathrm{ff}} \cdot \vec{n} &=  \vec{v}_{\mathrm{pm}} \cdot \vec{n}&&\text{on } \Sigma ,\label{eq:IC-mass} \\
    \label{eq:IC-momentum}
    -\vec{n} \cdot \ten{T}\left(\vec{v}_{\mathrm{ff}}, p_{\mathrm{ff}}\right) \cdot \vec{n} &= p_\mathrm{pm} &&\text{on } \Sigma,\\
    \label{eq:IC-BJ}
    \left(\vec{v}_{\mathrm{ff}}-\textbf{v}_{\mathrm{pm}}\right) \cdot \bm{\tau}-\alpha^{-1}\sqrt{K}\left(\left(\nabla \vec{v}_{\mathrm{ff}} + (\nabla \vec{v}_{\mathrm{ff}})^\top \right) \cdot \vec{n}\right) \cdot \bm{\tau} &=0 &&\text{on }\Sigma.
\end{align}
Here, $\vec{n}$ is the unit vector normal to the fluid--porous interface~$\Sigma$ pointing outward from the porous-medium domain $\Omega_{\mathrm{pm}}$, $\bm{\tau}$ is the unit vector tangential to the interface (Fig.~\ref{fig:domain-sparsity}), and $\alpha>0$ is the Beavers--Joseph slip coefficient. Different approaches to compute $\sqrt{K}$ exist in the literature. In this paper, we consider $\sqrt{K}= \sqrt{\bm{\tau}\cdot \ten{K}\cdot\bm{\tau}}$ as in~\cite{Rybak_etal_15}.

Saffman~\cite{Saffman} proposed a simplification of the Beavers--Joseph  condition~\eqref{eq:IC-BJ}, where the tangential porous-medium velocity $\textbf{v}_{\mathrm{pm}} \cdot \bm{\tau}$ on the interface is neglected
\begin{equation}
    \label{eq:IC-BJS}
    \textbf{v}_{\mathrm{ff}} \cdot \bm{\tau}-\alpha^{-1}\sqrt{K}\left(\left(\nabla\textbf{v}_{\mathrm{ff}} + (\nabla \textbf{v}_{\mathrm{ff}})^\top\right) \cdot \mathbf{n}\right) \cdot \bm{\tau} =0 \qquad \text{on }\Sigma.
\end{equation}
In the literature, the Stokes--Darcy problem with the Beavers--Joseph--Saffman condition~\eqref{eq:IC-BJS} is usually studied, both from the analytical and numerical point of view. Only a few papers focus on the original Beavers--Joseph condition~\eqref{eq:IC-BJ}. In this work, we develop and analyse preconditioners for the Stokes--Darcy problem \eqref{eq:Stokes}--\eqref{eq:BC-pm} with both sets of interface conditions, \eqref{eq:IC-mass}--\eqref{eq:IC-BJ} and \eqref{eq:IC-mass}, \eqref{eq:IC-momentum},~\eqref{eq:IC-BJS}.

\subsection{Discretisation}
\begin{figure}[!ht]
    \centering
    \includegraphics[width=0.6\linewidth]{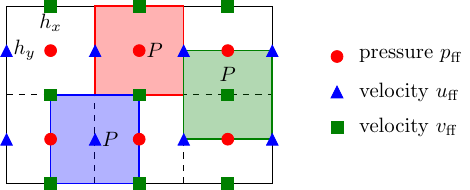}
    \caption{Staggered grid for the Stokes problem }
    \label{fig:staggered-grid-ff}
\end{figure}

The coupled Stokes--Darcy problems~\eqref{eq:Stokes}--\eqref{eq:IC-BJ} and \eqref{eq:Stokes}--\eqref{eq:IC-momentum}, \eqref{eq:IC-BJS} are discretised with the second order finite volume method. The MAC scheme (finite volume method on staggered grids) is used (Fig.~\ref{fig:staggered-grid-ff}) for the Stokes equations, e.g.~\cite{Eymard_etal_10,Rybak_etal_15,Greif_He_2023,Shiu-Ong-Lai-18,Luo_etal_17}. The porous-medium model~\eqref{eq:Darcy-law} is discretised in its primal form, where Darcy's law is substituted to the mass balance equation. Here, the pressure $p_\PM$ is the primary variable   which is defined in the control volume centres as well as on the fluid--porous interface and the external boundary of the domain. The porous-medium velocity components $\textbf{v}_\PM = (u_\PM;v_\PM)$ are computed on the control volume faces in a post-processing step.
This leads to the system of linear equations
\begin{equation}
\label{eq:system}
\mathcal{A}\textbf{x} = \textbf{b}, \qquad
\textbf{x}= (\textbf{v}_{\mathrm{ff}}; \, p_{\mathrm{ff}};\, p_{\mathrm{pm}})^\top, \qquad \mathcal{A} \in \{\mathcal{A}_{\mathrm{BJS}},\, \mathcal{A}_{\mathrm{BJ}}\},
\end{equation}
where $\textbf{v}_{\mathrm{ff}} \in \mathbb{R}^n, \, p_{\mathrm{ff}} \in \mathbb{R}^m,\, p_{\mathrm{pm}} \in \mathbb{R}^l$ are the primary variables, and the matrices are given by
\begin{equation}
 \label{eq:A-BJS-BJ}
     \mathcal{A}_{\mathrm{BJS}} = \left(\begin{array}{ccc}
     A & B^\top  & C^\top  \\
     B & 0 & 0 \\
     C & 0 & -D
 \end{array}\right), \qquad \mathcal{A}_{\mathrm{BJ}} = \left(\begin{array}{ccc}
     A & B^\top  & C_2^\top  \\
     B & 0 & 0 \\
     C_1 & 0 & -D
 \end{array}\right) .
 \end{equation}
Here, the blocks $A \in \mathbb{R}^{n\times n}$ and $D \in \mathbb{R}^{l\times l}$ are both symmetric and positive definite ($A = A^\top \succ 0, \, D = D^\top \succ 0$) and the matrix $B \in \mathbb{R}^{m \times n}$ has full row rank ($\text{rank}(B)~=~m)$. 
For the case of the Beavers--Joseph--Saffman condition~\eqref{eq:IC-BJS}, the matrix $\mathcal{A}_{\mathrm{BJS}}$ is symmetric ($\mathcal{A}_{\mathrm{BJS}} = \mathcal{A}_{\mathrm{BJS}}^{\top}$), however for the Beavers--Joseph condition~\eqref{eq:IC-BJ} it is not possible to get a completely symmetric matrix ($\mathcal{A}_{\mathrm{BJ}} \ne \mathcal{A}_{\mathrm{BJ}}^{\top}$ due to $C_2 = C_1+ E$, $E\neq 0$). Note that the matrix $\mathcal{A}_{\mathrm{BJS}}$ is a double saddle point matrix. The sparsity structure of $\mathcal{A}_{\mathrm{BJS}}$ and $\mathcal{A}_{\mathrm{BJ}}$ is presented in Fig.~\ref{fig:domain-sparsity} for the grid width $h_x=h_y=h=1/8$.

\begin{figure}[!ht]
    \centering
 \includegraphics[width=0.9\linewidth]{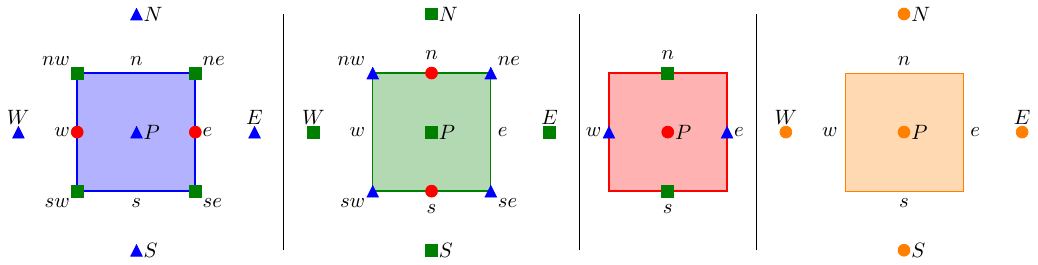}
    \caption{Control volumes for the primary variables $u_\FF$, $v_\FF$, $p_\FF$ and $p_\FF$}
    \label{fig:control-volume}
\end{figure}

To describe how the matrix blocks in~\eqref{eq:A-BJS-BJ} are built, we briefly provide the discretisation of the coupled Stokes--Darcy problem with emphasis on the interface. Note that in comparison to similar MAC schemes~\cite{Greif_He_2023, Shiu-Ong-Lai-18}, we have degrees of freedom for the free-flow velocities and the porous-medium pressure on the external boundaries of the domain and on the fluid--porous interface.

The blocks $A$ and $B^\top$ in~\eqref{eq:A-BJS-BJ} correspond to the discretised momentum balance equation in the Stokes system~\eqref{eq:Stokes}. The momentum equation in~\eqref{eq:Stokes} can be written for the horizontal and vertical component as
\begin{subequations}
    \begin{align}
    \label{eq:Stokes-MB-x}
        -\frac{\partial}{\partial x} \left(2 \mu \frac{\partial  u_\FF}{\partial x} \right) - \frac{\partial}{\partial y} \left(\mu \frac{\partial u_\FF}{\partial y} + \mu \frac{\partial v_\FF}{\partial x}  \right) + \frac{\partial p_\FF}{\partial x } &= f_{\FF}^u,\\
    \label{eq:Stokes-MB-y}
        - \frac{\partial}{\partial x} \left(\mu \frac{\partial v_\FF}{\partial x}  + \mu \frac{\partial u_\FF}{\partial y}  \right)-\frac{\partial}{\partial y} \left(2 \mu \frac{\partial v_\FF}{\partial y} \right) + \frac{\partial p_\FF}{\partial y } &= f_{\FF}^v .
    \end{align}
\end{subequations}

For the sake of compactness, we omit the subscripts ($\FF$, $\PM$) for the discretisation inside $\Omega_\FF$ and $\Omega_\PM$, and use them only on the interface $\Sigma$. For the inner control volume (Fig.~\ref{fig:control-volume}, blue) we get for the horizontal component of the momentum balance equation~\eqref{eq:Stokes-MB-x}:
\begin{equation}\label{eq:mom-u}
   (F_{x,e}^u -  F_{x,w}^u) \, h_y + (F_{y,n}^u - F_{y,s}^u) \, h_x = f^u h_x h_y,
\end{equation}
with the momentum fluxes
\begin{subequations}
\label{eq:fluxes}
\begin{align*}
    F_{x,e}^u = p_{e} - 2\mu \frac{u_{E}-u_{P}}{h_x}, &\quad F_{x,w}^u = p_{w} - 2 \mu \frac{u_{P}-u_{W}}{h_x},\\
    F_{y,n}^u = -\mu\frac{u_{N}-u_{P}}{h_y} - \mu \frac{v_{ne}-v_{nw}}{h_x}, &\quad F_{y,s}^u = -\mu \frac{u_{P}-u_{S}}{h_y}-\mu\frac{v_{se}-v_{sw}}{h_x}.
\end{align*}
\end{subequations}
Analogously, we have for the vertical component (Fig.~\ref{fig:control-volume}, green) of the momentum equation~\eqref{eq:Stokes-MB-y}:
\begin{equation}
   \label{eq:mom-v}
(F_{x,e}^v -  F_{x,w}^v) h_y + (F_{y,n}^v - F_{y,s}^v) h_x = f^v h_x h_y,
\end{equation}
with the momentum fluxes
\begin{subequations}
\begin{align*}
    F_{x,e}^v = -\mu\frac{u_{ne}-u_{se}}{h_y} - \mu \frac{v_{E}-v_{P}}{h_x}, &\quad F_{x,w}^v = -\mu \frac{u_{nw}-u_{sw}}{h_y}-\mu\frac{v_{P}-v_{W}}{h_x},\\
    F_{y,n}^v = p_{n} - 2\mu \frac{v_{N}-v_{P}}{h_y},&\quad F_{y,s}^v = p_{s} - 2 \mu \frac{v_{P}-v_{S}}{h_y}.
\end{align*}
\end{subequations}
The coefficients from the discretisation~\eqref{eq:mom-u} and~\eqref{eq:mom-v} appear in the blocks $A$ and~$B^\top$ in \eqref{eq:A-BJS-BJ}. 
The second row in~\eqref{eq:A-BJS-BJ} is the discrete form of the incompressibility condition (Fig.~\ref{fig:control-volume}, red):
\begin{equation}
    \label{eq:D_ff}
    \left(u_{e} - u_{w}\right)\, h_y+\left(v_{n} -v_{s}\right)\, h_x = 0.
\end{equation}

The third row is the discrete version of the porous-medium model~\eqref{eq:Darcy-law} in its primal form
\begin{align*}
-\nabla \cdot \left(\mu^{-1} \ten{K} \nabla p_\PM\right) = f_\PM.
\end{align*}
Here we label the discrete porous medium velocity and pressure with $\sim$. The scheme (Fig.~\ref{fig:control-volume}, orange) reads
\begin{equation}
    \label{eq:A_pm}
    (\tilde{u}_e-\tilde{u}_w)\, h_y + (\tilde{v}_n-\tilde{v}_s)\, h_x = f h_x h_y,
\end{equation}
where
\begin{align*}
    \tilde{u}_e = -\frac{k_{xx}}{\mu}\frac{\tilde{p}_{E}-\tilde{p}_{P}}{h_x}, & \quad 
    \tilde{u}_w = -\frac{k_{xx}}{\mu}\frac{\tilde{p}_{P}-\tilde{p}_{W}}{h_x},\\
    \tilde{v}_n = -\frac{k_{yy}}{\mu}\frac{\tilde{p}_{N}-\tilde{p}_{P}}{h_y}, & \quad
    \tilde{v}_s = -\frac{k_{yy}}{\mu}\frac{\tilde{p}_{P}-\tilde{p}_{S}}{h_y}.
\end{align*}
The discretisation~\eqref{eq:A_pm} enters the block $D$ in~\eqref{eq:A-BJS-BJ}.

\begin{figure}[!ht]
    \centering
    \includegraphics[width = 1.0\linewidth]{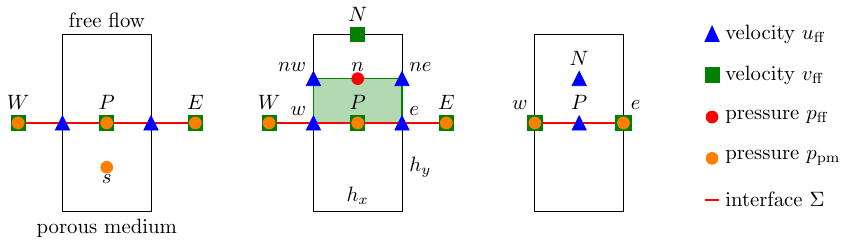}
    \caption{Stencils and primary variables for the coupling conditions}
    \label{fig:interface-stencil}
\end{figure}

Discretisation of coupling conditions is crucial for accuracy of the resulting coupled problem. To obtain the second order of convergence, we discretise the momentum equation near the interface using a half control volume (Fig.~\ref{fig:interface-stencil}, middle).   

The discretisation scheme (stencil in Fig.~\ref{fig:interface-stencil}, left) for the conservation of mass across the interface~\eqref{eq:IC-mass} reads
\begin{equation}
\label{eq:IC-mass-dis}
    -h_x\,v_{\FF,P} - 2\frac{k_{yy}}{\mu}\frac{h_x}{h_y} p_{\PM,s} + 2\frac{ k_{yy}}{\mu}\frac{h_x}{h_y} p_{\PM,P} = 0,
\end{equation}
where the coefficient $-h_x$ in the first term in \eqref{eq:IC-mass-dis} goes to the matrices $C$ and $C_1$, respectively, and the two other coefficients $\pm 2k_{yy}h_x/(\mu h_y)$ enter the matrix $D$ in~\eqref{eq:A-BJS-BJ}.

As mentioned above, we obtain the discrete form of the balance of normal forces \eqref{eq:IC-momentum} considering a control volume of size $h_x h_y/2$ (Fig.~\ref{fig:interface-stencil}, middle in green):
\begin{equation}
\label{eq:IC-BONF-dis}
\begin{aligned}
&\left(2\mu\frac{h_y}{h_x} + \mu\frac{h_x}{h_y}\right) v_{\FF,P}- 2\mu\frac{h_x}{h_y} v_{\FF,N}- \frac{\mu}{2} \frac{h_y}{h_x}v_{\FF,W} -\frac{\mu}{2}\frac{h_y}{h_x}v_{\FF,E} \\
&+\mu u_{\FF,nw}-\mu u_{\FF,ne}-\mu u_{\FF,w}+ \mu u_{\FF,e} + h_x p_{\FF,n}  -h_x p_{\PM,P} = f_P^v \frac{h_xh_y}{2}.
\end{aligned}
\end{equation}
The first eight coefficients in~\eqref{eq:IC-BONF-dis} contribute to the matrix $A$, the coefficient~$h_x$ to the matrix~$B$, and the coefficient $-h_x$ to the matrices $C^\top$ and $C_2^\top$, respectively. For discretisation of the Beavers--Joseph condition~\eqref{eq:IC-BJ} we get
\begin{equation}
\label{eq:IC-BJ-dis}
\begin{aligned}
    &\left(\mu\frac{\alpha}{\sqrt{k_{xx}}}h_x + 2\mu\frac{h_x}{h_y}\right) u_{\FF,P} - 2\mu \frac{h_x}{h_y}u_{\FF,N} +\mu v_{\FF,w} -\mu v_{\FF,e}\\
    &-\alpha\frac{k_{yy}}{\sqrt{k_{xx}}}p_{\PM,w} 
    + \alpha \frac{k_{yy}}{\sqrt{k_{xx}}}p_{\PM,e}  = 0,
\end{aligned}
\end{equation}
and for the simplification by Saffman~\eqref{eq:IC-BJS} we end up with
\begin{align}
\label{eq:IC-BJS-dis}
\left(\mu\frac{\alpha}{\sqrt{k_{xx}}}h_x + 2\mu\frac{h_x}{h_y}\right) u_{\FF,P} - 2\mu \frac{h_x}{h_y}u_{\FF,N} + \mu v_{\FF,w} - \mu v_{\FF,e} = 0.
\end{align}
The first four coefficients in~\eqref{eq:IC-BJ-dis} and~\eqref{eq:IC-BJS-dis} enter the matrix $A$. The fifth and sixth coefficients $\pm \alpha k_{yy}/\sqrt{k_{xx}}$ in~\eqref{eq:IC-BJ-dis} contribute to the matrix $E^\top$ and therefore to $C_2^\top$ breaking the symmetry of the matrix $\mathcal{A}_{\mathrm{BJ}}$. 

\section{Preconditioners}
 \label{sec:preconditioners}
 In this section, we propose three different types of preconditioners for coupled Stokes--Darcy problems, a block diagonal, a block triangular and a constraint one (section~\ref{sec:exact-preconditioners}). For efficient numerical simulations of large systems we develop inexact variants of these preconditioners (section~\ref{sec:inexact-preconditioners}).
 
 \subsection{Exact preconditioners}
 \label{sec:exact-preconditioners}
We solve the coupled flow problems \eqref{eq:system},~\eqref{eq:A-BJS-BJ} monolithically using the right-preconditioned GMRES method, which is applicable to non-symmetric matrices
\begin{equation}
\label{eq:right-preconditioning}
    \mathcal{A}\mathcal{P}^{-1} \overline{\textbf{x}} = \textbf{b},\qquad \overline{\textbf{x}}= \mathcal{P}\textbf{x}.
\end{equation}
As a first step in the construction of preconditioners, we decouple the Stokes--Darcy system and consider the following matrix
\begin{equation}
\label{eq:A-neglect-IC}
    \overline{\mathcal{A}} = \left(\begin{array}{ccc}
        A & B^\top  & 0\\
        B & 0 & 0 \\
        0 & 0 & -D
    \end{array}\right).
\end{equation} 
We develop preconditioners for the matrix $\overline{\mathcal{A}}$ and show theoretically and numerically that they are also suitable for the matrix $\mathcal{A} \in \{\mathcal{A}_{\mathrm{BJS}},\,\mathcal{A}_{\mathrm{BJ}}\}$ defined in~\eqref{eq:A-BJS-BJ}.
We propose block diagonal and block triangular preconditioners
\begin{equation}
\label{eq:P-D-P-T}
    \Pd = \left(\begin{array}{ccc}
   A  &  0 & 0 \\
   0  &  -S_B & 0 \\
   0 & 0 & -D
\end{array}\right), \qquad \Pt = \left( \begin{array}{ccc}
    A & B^\top  & 0 \\
    0 & -S_B & 0\\
    0 & 0 & -D
\end{array}\right),
\end{equation}
where $S_B:= BA^{-1}B^\top \in \mathbb{R}^{m \times m}$ is the Schur complement ($S_B = S_B^\top \succ 0$). A similar attempt to $\Pd$ was proposed in~\cite{Cai_etal_09}, where the Schur complement is approximated as $S_B \approx \mu^{-1} I$, where $I \in \mathbb{R}^{m\times m}$ is the identity matrix. Note that the block triangular preconditioner $\Pt$ given in~\eqref{eq:P-D-P-T} is a simplified version of the block triangular preconditioners proposed in~\cite{Beik_Benzi_18} for double saddle point problems. However, it justified its excellent performance in numerical simulations (section~\ref{sec:numerical-results}).

Furthermore, we construct a constraint preconditioner $\Pc$. 
The matrix $\overline{\mathcal{A}}$ can be interpreted as a standard saddle point matrix
\begin{equation}
\label{eq:A-bar}
    \overline{\mathcal{A}}=  \left(\begin{array}{cc}
   A  & \overline{B}^\top  \\
   \overline{B}  & - \overline{D}
\end{array}\right)\quad \text{with} \quad \overline{B} = \left(\begin{array}{c}
    B \\
     0 
\end{array}\right), \quad \overline{D} = \left(\begin{array}{cc}
   0  &  0\\
   0  &  D
\end{array}\right).
\end{equation}
Therefore, also preconditioners for standard saddle point problems  are suitable. We modify and generalise the following constraint preconditioner proposed in~\cite{Axelsson_Neytcheva_03}:
\begin{equation}
\label{eq:P_con}
    \Pc = \left(\begin{array}{cc}
   G  & \overline{B}^\top  \\
   \overline{B}  & -\overline{D}
\end{array}\right) \quad \text{with} \quad G = \text{diag}(A_{11},\, A_{22}), \quad G=G^\top \succ 0,
\end{equation}
where $G$ is a preconditioner to the block $A= \left(A_{ij}\right)_{i,j=1,2}$.

\subsection{Inexact preconditioners}
\label{sec:inexact-preconditioners}
To obtain accurate numerical results for the coupled Stokes--Darcy system with suitable boundary and interface conditions, we need to consider small grid widths which yield large linear systems. Since exact versions of preconditioners are computationally expensive, they have to be replaced by efficient inexact variants. The approximations of the blocks $A$, $G$, $S_B$ and $D$ in~\eqref{eq:P-D-P-T} and~\eqref{eq:P_con} are marked by $\hat{A},\, \hat{G},\, \hat{S}_B$ and $\hat{D}$, respectively, and the inexact versions of the corresponding preconditioners are $\hPd,\, \hPt $ and $\hPc$. The approximations should be easily invertible to reduce the computational effort. 
We replace the inverses of the blocks $A$, $G$ and $D$ from~\eqref{eq:P-D-P-T} and~\eqref{eq:P_con} applied in~\eqref{eq:right-preconditioning} with the approximations $\hat{A}^{-1} = \hat{G}^{-1} = \text{diag}\left(\text{AMG}(A_{11}) ,\text{AMG}(A_{22})\right)$ and $\hat{D}^{-1} = \text{AMG}(D)$, respectively. Here, $\text{AMG}(A_{11})$, $\text{AMG}(A_{22})$ and $\text{AMG}(D)$ are algebraic multigrid methods. 

For the Stokes--Darcy system discretised with FEM a common way to approximate the Schur complement $S_B$ is to use the Stokes pressure mass matrix, e.g.~\cite[Theorem~3.22]{Finite-Elements-and-Fast-Iterative-Solvers}.
To define the finite volume analogon of the pressure mass matrix, we consider the continuous version of the Schur complement
\begin{equation}
    \mathcal{S}_B := -\nabla \cdot \left[-\nabla \cdot \left(\mu \left(\nabla + \nabla^\top \right)\right)\right]^{-1} \nabla.
\end{equation}
We note that $\mathcal{S}_B \approx (2\mu)^{-1} I$. 
Since we use the finite volume scheme, we approximate the Schur complement $S_B$ by $\hat{S}_B = (2\mu)^{-1}h_xh_y I$.

\section{Spectral analysis}
\label{sec:spectral-analysis}
 The goal of this section is to analyse the spectra of the preconditioned matrices $\mathcal{A}\Pd^{-1}, \, \mathcal{A}\Pt^{-1}$ and $\mathcal{A}\Pc^{-1}$ for $\mathcal{A} \in \{\mathcal{A}_{\mathrm{BJS}}, \mathcal{A}_{\mathrm{BJ}}\}$ defined in~\eqref{eq:A-BJS-BJ}. Eigenvalues clustered around one and/or a clustered spectrum away from zero often provide fast convergence of the Krylov subspace methods~\cite{Benzi_02}. 

For the matrix $H = H^\top \succ 0$, we define the vector norm and the corresponding induced matrix norm
 $$\langle x, x \rangle_H = \langle Hx, x\rangle = x^\top H x = \|x\|^2_H, \qquad \|M\|_{H} = \max_{x \neq 0} \frac{\|Mx\|_{H}}{\|x\|_{H}}.$$
 The extension of the induced matrix norm for $H_1 = H_1^\top \succ 0,\, H_2 = H_2^\top \succ 0$ is given by~\cite[Problem 5.6.P4]{Matrix-analysis}:
 \begin{equation}
 \label{eq:ext-mat-norm}
     \|M\|_{H_1,H_2} = \max_{x \in \mathbb{R}^{n}\backslash\{0\}} \frac{\|Mx\|_{H_2}}{\|x\|_{H_1}}, \qquad M \in \mathbb{R}^{m \times n}.
 \end{equation}
 Moreover, we have the following equalities
\begin{equation}
\label{eq:norm-equalities}
 \!\!  \|H_2^{-1/2} M H_1^{-1/2}\|_2 = \|M\|_{H_1,H_2^{-1}} = \|M H_1^{-1}\|_{H_1^{-1},H_2^{-1}} = \|H_2^{-1}M\|_{H_1,H_2}.  \! \!   
\end{equation}
 Note that for $H_1 = H_2 = H$ and $m = n$, we get the standard induced matrix norm $\|M\|_{H,H} = \|M\|_H$.

In the following, we prove three practical theorems (Theorems~\ref{theo:eig-APd-2},~\ref{theo:eig-APt-2},~\ref{theo:eig-APc}) on the clustering of the eigenvalues for the preconditioned matrices $\ABJS \mathcal{P}^{-1}$ and $\ABJ \mathcal{P}^{-1}$ with $\mathcal{P} \in \{\Pd, \Pt, \Pc\}$. Furthermore, we provide bounds on the spectra of $\ABJS \Pd^{-1}$ and $\ABJS \Pt^{-1}$ (Theorems~\ref{theo:eig-APd},~\ref{theo:eig-APt}). Note that the eigenvalues of the matrices $\ABJS \Pc^{-1}$ and $\ABJ \Pc^{-1}$ are already determined in Theorem~\ref{theo:eig-APc}. We start by proving a supplementary lemma needed in Theorems~\ref{theo:eig-APd} and~\ref{theo:eig-APt}.

\begin{lemma}
\label{lem:est-rq}
Let $A = A^\top \succ 0$ and $B$ has full rank, then the following equalities
\begin{enumerate}[(i)]
    \item $\displaystyle \max_{x \neq 0} \frac{x^\top B^\top (BA^{-1}B^\top)^{-1}B x}{x^\top A x} = 1 $,
    \item $\displaystyle \frac{x^\top B^\top (BA^{-1}B^\top)^{-1}Bx}{x^\top A x} = 1 \quad \text{for } x \not\in \operatorname{ker}(BA^{-1})$
\end{enumerate}
hold.
 \end{lemma}
 \begin{proof}
  \begin{enumerate}[$(i)$]
    \item Simple algebraic manipulations yield
 \begin{align*}
        &\max_{x \neq 0} \frac{x^\top B^\top S_B^{-1}B x}{x^\top A x} \overset{~\eqref{eq:ext-mat-norm}}{=}\max_{x \neq 0} \frac{\|Bx\|^2_{S_B^{-1}}}{\|x\|^2_{A}}
        = \|B\|^2_{A, S_B^{-1}} \\
     &\overset{\eqref{eq:norm-equalities}}{=} \left\|S_B^{-1/2}BA^{-1/2}\right\|^2_2 
        \!= \! \left\|A^{-1/2}B^\top S_B^{-1/2}\right\|^2_2 \\
        &\overset{\eqref{eq:norm-equalities}}{=} \left\|B^\top\right\|^2_{S_B,A^{-1}} \!=\! 
         \max_{y \neq 0} \frac{\|B^\top y\|^2_{A^{-1}}}{\|y\|^2_{S_B}}\overset{\eqref{eq:ext-mat-norm}}{=}\max_{y \neq 0} \frac{y^\top S_B y}{y^\top S_B y} = 1.
\end{align*}
\item We obtain
     $$\frac{\langle B^\top S_B^{-1} B x, x\rangle}{\langle Ax, x\rangle} = \frac{\langle B^\top S_B^{-1} B A^{-1}y, y\rangle}{\langle y,y\rangle}.$$
     We have $B^\top S_B^{-1} BA^{-1} = I$ on $\text{range}(B^\top) \supset \text{range}(A^{-1}B^\top)$, see~\cite[Section 1]{Elman_99}. Thus  due to $\mathbb{R}^m = \text{range}(A^{-1}B^\top) \oplus \text{ker}(BA^{-1})$, this equality also holds for $\mathbb{R}^m\backslash\operatorname{ker}(BA^{-1})$.
\end{enumerate}
\hfill~\qed
\end{proof}
In the following, we define a column vector $(x;y;z)^\top \in \mathbb{C}^{n+m+l}$ for given vectors $x \in \mathbb{C}^n,\, y \in \mathbb{C}^m$ and $z \in \mathbb{C}^l$.
\begin{theorem}
 \label{theo:eig-APd-2}
The eigenvalues of the preconditioned matrix $\ABJS \Pd^{-1}$ with the block diagonal preconditioner $\Pd$ defined in~\eqref{eq:P-D-P-T} are $\lambda = 1$ or cluster around $(1 \pm \sqrt{3}i)/2$ and $1$ for $h \rightarrow 0$. The same clustering holds for $\ABJ \Pd^{-1}$ under the additional assumption $\alpha k_{yy} / \sqrt{k_{xx}} \ll 1$.
 \end{theorem}
 \begin{proof}
 We carry out the proof for the general case $\ABJ \Pd^{-1}$.
  Let $\lambda$ be an eigenvalue of $\mathcal{A}_{\mathrm{BJ}}\Pd^{-1}$ to the eigenvector $(x;y;z)^\top \neq 0$ such that
    \begin{align}
        Ax + B^\top y + C_2^\top z &= \lambda Ax, \label{eq:eig-APd-1-2}\\
        Bx &= -\lambda S_B y, \label{eq:eig-APd-2-2}\\
        C_1 x - D z &= -\lambda Dz. \label{eq:eig-APd-3-2}
    \end{align}
    We obtain $\lambda = 1$ with the corresponding eigenvectors $(x;0;0)^\top$ for $0 \neq x \in \text{ker}(B) \cap \text{ker}(C_1)$ and $(0;0;z)^\top$ for $0 \neq z \in \text{ker}(C_2^\top)$.

   Now we assume $\lambda \neq 1$. To show clustering around $(1 \pm \sqrt{3}i)/2$ and $1$, we get \mbox{$x =(A^{-1}B^\top y + A^{-1} C_2^\top z)/(\lambda-1)$} from~\eqref{eq:eig-APd-1-2} and insert it into~\eqref{eq:eig-APd-2-2}. Rearranging the terms and multiplying the expression with $y^\top$ from the left, we obtain
    \begin{equation*}
     a \lambda^2 - a\lambda + a+(b + e) = 0,
     \end{equation*}   
     \begin{equation}   
     \!\!   a = y^\top S_B y > 0, \;
        b = y^\top BA^{-1}C_1^\top z = z^\top C_1A^{-1}B^\top y,\; e = y^\top BA^{-1}E^\top z. \!\!\label{eq:a-b}
     \end{equation}
    The solution of this equation is
    \begin{equation}
    \label{eq:clust-APd-1}
      \lambda_{1,2} = \left(1 \pm \sqrt{-3 -4(b+e)/a}\right)/2. 
    \end{equation}
    Analogously, inserting \mbox{$x =(A^{-1}B^\top y + A^{-1} C_2^\top z)/(\lambda-1)$} in~\eqref{eq:eig-APd-3-2}, rearranging the terms and multiplying from the left with $z^\top$, we get
    \begin{equation}
    \label{eq:c-d}
        c \lambda^2 - 2c \lambda +c+b+d = 0, \quad c = z^\top D z > 0, \; d = z^\top C_1A^{-1}C_2^\top z.
    \end{equation}
    Solving the equation above leads to
    \begin{equation}
    \label{eq:clust-APd-2}
        \lambda_{1,2} = 1\pm \sqrt{-(b+d)/c}.
    \end{equation}
    The terms $b$ and $d$ defined in~\eqref{eq:a-b} and~\eqref{eq:c-d} converge to zero for $h\rightarrow 0$, because the entries of $C_1$ are of order $O(h)$ due to discretisation of the interface conditions~\eqref{eq:IC-mass-dis} and~\eqref{eq:IC-BONF-dis}. Due to the assumption $\alpha k_{yy}/\sqrt{k_{xx}} \ll 1$, it holds $e \ll 1$ (see discretisation in~\eqref{eq:IC-BJ-dis}). Therefore, the eigenvalues cluster around $(1 \pm \sqrt{3}i)/2$ and $1$. 
    Note that in the case of $\ABJS \Pd^{-1}$, we have $C_1 = C_2 = C$ and therefore $e=0$.
    \hfill~\qed
 \end{proof}
 The eigenvalue distributions of the preconditioned matrices $\ABJS\Pd^{-1}$ and $\ABJ\Pd^{-1}$ provided in Fig.~\ref{fig:eigenvalues} confirm the theoretical results obtained in Theorem~\ref{theo:eig-APd-2}.
  
 \begin{theorem}
 \label{theo:eig-APd}
The preconditioned matrix $\ABJS \Pd^{-1}$ with the block diagonal preconditioner $\Pd$ defined in~\eqref{eq:P-D-P-T} has either the eigenvalue $\lambda = 1$ or $|\lambda|\geq \tau$ and $\zeta \leq |\lambda -1| \leq 1 + \omega$ for $\tau, \, \zeta,\, \omega > 0$.
 \end{theorem}
 \begin{proof}
  Let $\lambda$ be an eigenvalue of $\mathcal{A}_{\mathrm{BJS}}\Pd^{-1}$ to the eigenvector $(x;y;z)^\top \neq 0$:
    \begin{align}
        Ax + B^\top y + C^\top z &= \lambda Ax, \label{eq:eig-APd-1}\\
        Bx &= -\lambda S_B y, \label{eq:eig-APd-2}\\
        C x - D z &= -\lambda Dz. \label{eq:eig-APd-3}
    \end{align}
   Considering eigenvectors $(x;0;0)^\top$ for $0 \neq x \in \text{ker}(B) \cap \text{ker}(C)$ and $(0;0;z)^\top$ for $0 \neq z \in \text{ker}(C^\top)$, we obtain $\lambda = 1$.

   Now we assume $\lambda \neq 1$. From~\eqref{eq:eig-APd-2} and~\eqref{eq:eig-APd-3}, we get $y=-\lambda^{-1}S_B^{-1} B x$ and $z = 1/(1-\lambda)D^{-1}C x$. Substitution of these vectors in~\eqref{eq:eig-APd-1} yields 
    \begin{equation}
        \label{eq:char-APd}
        \lambda (1-\lambda)^2 q - (1-\lambda) r + \lambda p = 0,
    \end{equation}
    where
    \begin{equation}
    \label{eq:q-p-r}
        q = x^\top A x > 0, \quad r=x^\top B^\top S_B^{-1} B x \geq 0, \quad p=x^\top C^\top D^{-1}C x \geq 0.
    \end{equation}
    For $p= 0$ ($x \in \text{ker}(C)$), we obtain $\lambda = 1/2\pm \sqrt{1/4-r/q}$.\\
    For $r=0$ $(x\in \text{ker}(B))$, we get $\lambda = 1 \pm \sqrt{-p/q}$.\\
    For $r \neq 0,\, p \neq 0$, substitution of $\lambda = t+1$ in~\eqref{eq:char-APd} yields
    \begin{equation}
    \label{eq:char-APd-t}
        qt^3 + qt^2+(r+p)t+p = 0.
    \end{equation}
    Applying~\cite[Theorem 2.5]{Beik_Benzi_22}, we get the bounds
    $$ 0<\zeta:=\min\left\{\frac{p}{p+r},\,\frac{r+p}{q},\, 1\right\} \leq |\lambda -1| \leq \max\left\{\frac{p}{p+r},\,\frac{r+p}{q},\, 1\right\}.$$
    For the upper bound, we show that
    $$\max\left\{\frac{p}{p+r},\,\frac{r+p}{q},\, 1\right\} \leq \max\left\{\max_{x\neq 0} \frac{p}{p+r},\, \max_{x\neq 0} \frac{r+p}{q},\, 1\right\} \leq 1+\omega.$$
    Here we applied Lemma~\ref{lem:est-rq}$(i)$ with $r$ and $q$ defined in~\eqref{eq:q-p-r} to obtain  \mbox{$\displaystyle \max\limits_{x\neq 0}(r/q) = 1$}, leading to 
    \[\max\left\{\max_{x\neq 0} \frac{p}{p+r},\, \max_{x\neq 0} \frac{r+p}{q},\, 1\right\} \leq 1 + \max_{x \neq 0} \frac{p}{q} =: 1 + \omega.\]
    Using~\eqref{eq:char-APd} we estimate
    \begin{align*}
        |\lambda| = |t+1| \geq \frac{r|t|}{q|t|^2 + p} \geq \frac{r \zeta}{q(1+\omega)^2 +p} =: \tau,
    \end{align*}
    which gives us the lower bound on $|\lambda|$. \hfill~\qed
 \end{proof}

\begin{theorem}
\label{theo:eig-APt-2}
The eigenvalues of the preconditioned matrix $\ABJS \Pt^{-1}$ with the block triangular preconditioner $\Pt$ defined in~\eqref{eq:P-D-P-T} are $\lambda = 1$ or cluster around $1$ for $h \rightarrow 0$. The same clustering holds for $\ABJ \Pt^{-1}$ under the additional assumption $\alpha k_{yy} / \sqrt{k_{xx}} \ll 1$.
\end{theorem}
\begin{proof}
    Let $\lambda$ be an eigenvalue of $\mathcal{A}_{\mathrm{BJ}}\Pt^{-1}$ to the eigenvector $(x;y;z)^\top \neq 0$ such that
    \begin{align}
        Ax + B^\top y + C_2^\top z &= \lambda (Ax + B^\top y), \label{eq:eig-APt-1-2}\\
        Bx &= -\lambda S_B y, \label{eq:eig-APt-2-2}\\
        C_1 x - D z &= -\lambda Dz. \label{eq:eig-APt-3-2}
    \end{align}
    We get $\lambda = 1$ with the corresponding eigenvectors $(x;0;0)^\top$ for $0 \neq x \in \text{ker}(B) \cap \text{ker}(C_1)$ and $(0;0;z)^\top$ for $0 \neq z \in \text{ker}(C_2^\top)$. Now, we assume $\lambda \neq 1$. 
    To show clustering around $1$, we get $x=-A^{-1}B^\top y - A^{-1}C_2^\top z/(1-\lambda)$ from~\eqref{eq:eig-APt-1-2}, insert it into~\eqref{eq:eig-APt-2-2} and~\eqref{eq:eig-APt-3-2} and multiply the resulting expressions with $y^\top$ and $z^\top$, respectively. This leads to
    \begin{align}
    \label{eq:lambda-Pt-1}
        a \lambda^2 - 2a \lambda +a+(b+e) = 0,\\
        a= y^\top S_B y,\; b = y^\top BA^{-1}C_1^\top z  = z^\top C_1A^{-1}B^\top y,\; e = y^\top BA^{-1}E^\top z
    \end{align}
    and
    \begin{align}
     \label{eq:lambda-Pt-2}
        c \lambda^2 - (2c+b) \lambda +c + b +d = 0, \qquad c= z^\top D z,\quad d= z^\top C_1A^{-1}C_2^\top z.
    \end{align}
    Solutions of equations~\eqref{eq:lambda-Pt-1} and~\eqref{eq:lambda-Pt-2} are
    \begin{align*}
        \lambda_{1,2} = 1\pm \sqrt{-(b+e)/a},\qquad  \lambda_{1,2} = 1 + \left(b/c \pm \sqrt{(b/c)^2-4d/c}\right)/2.
    \end{align*}
    As in Theorem~\ref{theo:eig-APd-2}, we get $e \ll 1$ under the assumption $\alpha k_{yy}/\sqrt{k_{xx}}$. Since $b$ and $d$ converge to zero for $h \rightarrow 0$, the eigenvalues are clustered around $1$. 
    Note that in the case of $\ABJS \Pt^{-1}$, we have $C_1 = C_2 = C$ and therefore $e=0$. \hfill~\qed
\end{proof}
In Fig.~\ref{fig:eigenvalues}, we observe clustering around $1$ for both matrices $\ABJS \Pt^{-1}$ and $\ABJ \Pt^{-1}$, that validates numerically the theoretical results in Theorem~\ref{theo:eig-APt-2}.

\begin{theorem}
\label{theo:eig-APt}
    The preconditioned matrix $\ABJS \Pt^{-1}$ defined in~\eqref{eq:A-BJS-BJ} with the block triangular preconditioner $\Pt$ defined in~\eqref{eq:P-D-P-T} has either the eigenvalue $\lambda = 1$ or $|\lambda|\geq \tau$ and $\zeta \leq |\lambda -2| \leq 3$ for $\zeta, \, \tau > 0$.   
\end{theorem}
\begin{proof}
    Let $\lambda$ be an eigenvalue of $\mathcal{A}_{\mathrm{BJS}}\Pt^{-1}$ to the eigenvector $(x;y;z)^\top \neq 0$ such that
    \begin{align}
        Ax + B^\top y + C^\top z &= \lambda (Ax + B^\top y), \label{eq:eig-APt-1}\\
        Bx &= -\lambda S_B y, \label{eq:eig-APt-2}\\
        C x - D z &= -\lambda Dz. \label{eq:eig-APt-3}
    \end{align}
    We get $\lambda = 1$ with the corresponding eigenvectors $(x;0;0)^\top$ for $0 \neq x \in \text{ker}(B) \cap \text{ker}(C)$ and $(0;0;z)^\top$ for $0 \neq z \in \text{ker}(C^\top)$. Now, we assume $\lambda \neq 1$.
    From~\eqref{eq:eig-APt-2} and~\eqref{eq:eig-APt-3} we get $y=-S_B^{-1} B x / \lambda$ and $z = D^{-1}C x / (1-\lambda)$. Inserting $y$ and $z$ in~\eqref{eq:eig-APt-1} and rearranging the terms, we obtain
    \begin{equation}
        \label{eq:char-APt}
        \lambda(1-\lambda)^2 q  -  (1-\lambda)^2 r + \lambda p = 0,
    \end{equation}
    where $q,r$ and $p$ are defined in~\eqref{eq:q-p-r}.\\
    For $p = 0$ ($x \in \text{ker}(C)$), we get $\lambda = r/q$.\\
    For $r= 0$ ($x \in \text{ker}(B)$), we get $\lambda = 1 \pm \sqrt{-p/q}$.\\
    For $p \neq 0,\, r\neq 0$, substitution of $\lambda = t+2$ in~\eqref{eq:char-APt} yields
    \begin{equation}\label{eq:triag-quadr}
    qt^3 + (4q -r)t^2 + (5q-2r+p)t+2q-r+2p = 0.
    \end{equation}
     Following a similar procedure as in Theorem~\ref{theo:eig-APd} and taking Lemma~\ref{lem:est-rq}$(ii)$ into account, we obtain 
    \begin{equation}
        \zeta:= \min\left\{ \frac{2q - r+2p}{5q-2r+p}, \frac{5q-2r+p}{4q-r}\right\} \leq |\lambda - 2 |\leq \frac{4q - r}{q} = 3,
    \end{equation}
    since the coefficients in~\eqref{eq:triag-quadr} are positive and $q/r = 1$ for $r\neq 0$.
    We obtain the lower bound on $|\lambda|$ from~\eqref{eq:char-APt} with $\lambda= s+1$:
    \begin{align*}
        |\lambda| = |s+1| \geq \frac{r|s|^2}{q|s|^2 + p}  =: \tau. 
    \end{align*}
    \hfill~\qed 
\end{proof}

\begin{theorem}
\label{theo:eig-APc}
The preconditioned matrices $\ABJS \Pc^{-1}$ and $\ABJ \Pc^{-1}$ have the eigenvalues either $\lambda = 1$ or $\lambda = \left( 1 +\eta \pm \sqrt{(\eta - 1)^2-4\xi}\right)/2$ with $\xi = O(h)$ and $\eta = \left(x^\top Ax\right)/\left(x^\top G x\right) > 0$.
\end{theorem}
\begin{proof}
    Let $\lambda$ be an eigenvalue of $\ABJ\Pc^{-1}$ to the eigenvector $(x;y;z)^\top \neq 0$ such that
    \begin{align}
        Ax + B^\top y + C_2^\top z &= \lambda (Gx + B^\top y), \label{eq:eig-1}\\
        Bx &= \lambda Bx, \label{eq:eig-2}\\
        C_1 x - D z &= -\lambda Dz. \label{eq:eig-3}
    \end{align}
    We get $\lambda = 1$ with the corresponding eigenvector $(0;0;z)^\top$ for $0 \neq z \in \text{ker}(C_2^\top)$.
    If $\lambda \neq 1$, it holds $Bx = 0$. We obtain $z = D^{-1} C_1 x / (1-\lambda)$ from~\eqref{eq:eig-3}, substitute it into~\eqref{eq:eig-1} and multiply the resulting expression from the left with~$x^\top$. Rearrangement of the terms yields 
    \begin{equation}
    \label{eq:xi}
        \lambda^2 - \lambda(\eta + 1) + (\eta + \xi)= 0, \qquad \eta = \frac{x^\top A x}{x^\top G x} > 0, \quad \xi = \frac{x^\top C_2^\top D^{-1} C_1 x}{x^\top G x},
    \end{equation}
    with the roots
    \[\lambda_{1,2} = \left( 1 +\eta \pm \sqrt{(\eta - 1)^2-4\xi}\right)/2.\]
    Note that $\xi$ defined in~\eqref{eq:xi} is of order $O(h)$, because the entries of $C_1$ are of order $O(h)$. 
    \hfill~\qed
\end{proof}
In Fig.~\ref{fig:eigenvalues} and Tab.~\ref{tab:cluster-eigenvalues-Pc}, we observe clusters around $1$ and $\eta$, where $\eta$ depends on the choice of the block~$G$ in~\eqref{eq:P_con}. This is in accordance with the theoretical results of Theorem~\ref{theo:eig-APc}.

\section{Numerical study}
\label{sec:numerical-results}
In this section, we present numerical simulation results for two coupled Stokes--Darcy problems~\eqref{eq:Stokes}--\eqref{eq:BC-pm}: $(i)$~\emph{Problem $\ABJS$} is completed with the conservation of mass~\eqref{eq:IC-mass}, the balance of normal forces~\eqref{eq:IC-momentum} and the Beavers--Joseph--Saffman condition~\eqref{eq:IC-BJS} on the tangential velocity, and $(ii)$~\emph{Problem~$\ABJ$} with the coupling conditions~\eqref{eq:IC-mass},~\eqref{eq:IC-momentum} and the Beavers--Joseph condition~\eqref{eq:IC-BJ}.

The Stokes--Darcy problems 
are implemented using our in-house C++ code. To compute the eigenvalues, we use the \textsc{NumPy} function \code{linalg.eig} in Python 3.10. We solve the original and the right-preconditioned systems with the GMRES method considering both, the exact and inexact versions of the preconditioners. The AMG method for the inexact versions is implemented using the MATLAB toolbox IFISS~\cite{iffis}. 

The initial guess is always taken $\vec{x}_0 = \vec{0}$. The stopping criterion is either the maximum number of iteration steps $n_{\max}=2000$ or  $\|\mathcal{A}\vec{x}_n - \vec{b}\|_2 \leq \varepsilon_{\mathrm{tol}}$ for the tolerance $\varepsilon_{\mathrm{tol}} = 10^{-8}$. All computations are carried out on a laptop with an 12th Gen Intel(R) Core(TM) i7 1255U processor and 2$\times$16GB RAM using MATLAB.R2019b.

\subsection{Convergence order of the MAC scheme}

In this section, we study numerically the convergence order of the MAC scheme for the coupled Stokes--Darcy problems (\textit{Problem} $\ABJS$ and \textit{Problem} $\ABJ$). We have developed the second order method for the coupled problem using a different discretisation near the fluid--porous interface compared to~\cite{Greif_He_2023,Shiu-Ong-Lai-18}, where the overall convergence order is one. To demonstrate the advantage of our discretisation scheme, we consider the example (Test~1) from~\cite[Example~2]{Greif_He_2023} and~\cite[Example~2]{Shiu-Ong-Lai-18}, and provide some more additional tests.

\enlargethispage{0.25cm}

\begin{test}
\label{ex:1}
The flow domains are $\Omega_\FF = [0,1] \times [1,2]$ and $\Omega_\PM = [0,1] \times [0,1]$ with the fluid--porous interface $\Sigma = [0,1] \times\{1\}$. The porous medium is isotropic, $\ten K = k \ten I$. All physical parameters are $\alpha = \mu = k = 1$. The analytical solution is given by
\begin{align*}
    u^{\ex}_{\FF} &= (y-1)^2+x(y-1)+3x-1, \\
    v^{\ex}_{\FF} &= x(x-1)-0.5(y-1)^2-3y+1, \\
    p^{\ex}_{\FF} &= 2x+y-1,\\
    p^{\ex}_{\PM} &= x(1-x)(y-1)+(y-1)^3/3+2x+2y+4,
\end{align*}
yielding the right-hand sides
$\textbf{f}_{\FF} = \textbf{0}$, $f_{\PM} = 0$. The coupling condition on the tangential velocity is the Beavers--Joseph--Saffman version~\eqref{eq:IC-BJS} as in~\cite{Greif_He_2023,Shiu-Ong-Lai-18}, and on the external boundary the corresponding Dirichlet conditions are set.
\end{test}

In Tab.~\ref{tab:L2-error-test-1}, we provide the errors for all primary variables 
\[
\varepsilon_{w} = \|w - w^\ex\|_{L_2}, \quad w \in \{u_\FF, v_\FF, p_\FF, p_\PM\}.
\]
The corresponding convergence order is given in Tab.~\ref{tab:convergence-rates-ex1}. 
For the MAC scheme used in this work, the dimension of the linear systems~\eqref{eq:system},~\eqref{eq:A-BJS-BJ} is determined as 
\[d = (n_x+1) \cdot (n_y+2) + (n_x+2) \cdot (n_y+1) + n_x \cdot n_y + (m_x+2)\cdot (m_y+2),\]
where $n_x = n_y = 1/h$ is the number of control volumes in the free-flow domain and $m_x = m_y = 1/h$ is the number of control volumes in the porous medium in $x$- and $y$-direction, respectively.
Note that we have the second order of convergence for all primary variables, whereas in~\cite{Greif_He_2023,Shiu-Ong-Lai-18} only the first order is achieved for the same problem.

\begin{table}[!ht]
\caption{Errors for all primary variables in Test~\ref{ex:1}}
    \centering
    \begin{tabular}{|c||c|c|c|c|c|c|}
    \hline
       Grid & $h=1/8$ & $h=1/16$ & $h=1/32$ & $h=1/64$ & $h=1/128$ & $h=1/256$\\
       \hline
       Dim & $d=344$ & $d=1192$ & $d=4424$ & $d=17032$ & $d=66824$ & $d=264712$\\
       \hline \hline
       $\varepsilon_{u_\FF}$ & 
       9.3098e-4 & 2.3493e-4 & 5.9117e-5 & 1.4837e-5 & 3.7188e-6 & 9.3118e-7\\
       \hline
       $\varepsilon_{v_\FF}$ & 1.4285e-3 & 3.8177e-4 & 9.8864e-5 & 2.5182e-5 & 6.3565e-6 & 1.5943e-6\\
       \hline
       $\varepsilon_{p_\FF}$ & 3.2984e-2 & 9.4550e-3 & 2.6292e-3 & 7.1738e-4 & 1.9318e-4 & 5.1522e-5\\
       \hline 
       $\varepsilon_{p_\PM}$ & 1.1780e-3 & 3.2131e-4 & 8.3900e-5 & 2.1453e-5& 5.4261e-6 & 1.3647e-6\\
       \hline
    \end{tabular}
    \label{tab:L2-error-test-1}
\end{table}
\begin{table}[!ht]
\caption{Convergence order for all primary variables in Test~\ref{ex:1}}
    \centering
    \begin{tabular}{|c||c|c|c|c|c|}
    \hline
       Grid 1/Grid 2  & $8/16$ & $16/32$ & $32/64$ & $64/128$ & $128/256$ \\
       \hline\hline
    $u_\FF$ & 1.9865 & 1.9885 & 1.9905 &  1.9919 & 1.9931 \\
    \hline
    $v_\FF$ & 1.9038 & 1.9265 & 1.9420 & 1.9530 & 1.9610 \\
    \hline 
    $p_\FF$ & 1.8026 & 1.8245 & 1.8410 & 1.8539 & 1.8645 \\
    \hline
    $p_\PM$ & 1.8743 & 1.9058 & 1.9263 & 1.9406 & 1.9507 \\
    \hline
    \end{tabular}
    \label{tab:convergence-rates-ex1}
\end{table}
We developed an additional test, which is used to investigate both the convergence of the MAC scheme and the performance of the developed preconditioners with respect to different physical parameters. Here, both versions of the coupling condition on the tangential velocity  are considered: the Beavers--Joseph--Saffman (\textit{Problem} $\ABJS$) and the more general Beavers--Joseph condition (\textit{Problem} $\ABJ$).

\begin{test}
    \label{ex:3}
   The coupled flow domain is  $\Omega_\FF \cup \Sigma\cup \Omega_\PM$, where $\Omega_\FF = [0,1] \times [1,2]$, $\Omega_\PM = [0,1] \times [0,1]$ and $\Sigma = [0,1] \times\{1\}$. The porous medium is isotropic, $\ten K = k \ten I$, and there are no constraints on the physical parameters $\alpha$, $\mu$ and $k$. The analytical solution, which satisfies the coupling conditions, is given by
    \begin{align*}
        u^\ex_{\FF} &= -\text{cos}(\pi x) \text{sin}(\pi y), \quad   &&v^\ex_\FF = \text{sin}(\pi x) \text{cos}(\pi y),\\
        p^\ex_\FF &= \mu k^{-1} (y-1) \text{sin}(\pi x), \quad
        &&p^\ex_\PM = \mu k^{-1} (y^2 - y) \text{sin}(\pi x).
    \end{align*}
We consider Dirichlet boundary conditions on the external boundary and compute the right-hand sides substituting the exact solution in the corresponding equations. 
\end{test}
\enlargethispage*{1.7cm}
We provide the errors $\epsilon_w$ for all primary variables $w\in \{u_\FF, v_\FF, p_\FF,p_\PM\}$ for \textit{Problem} $\ABJS$ (Tab.~\ref{tab:error-2-BJS}) and \textit{Problem} $\ABJ$ (Tab.~\ref{tab:error-2-BJ}) with physical parameters $\alpha = 1$, $\mu = 10^{-3}$ and $k=10^{-2}$, and the corresponding convergence rates in Tab.~\ref{tab:convergence-rates-ex-2-BJS} and Tab.~\ref{tab:convergence-rates-ex-2-BJ}, accordingly. In both cases, the second order convergence is confirmed.
\begin{table}[!ht]
    \caption{Errors for all primary variables in Test~\ref{ex:3} (\textit{Problem} $\ABJS$) for physical parameters $\alpha = 1$, $\mu=10^{-3}$ and $k=10^{-2}$}
    \centering
    \begin{tabular}{|c||c|c|c|c|c|c|}
    \hline
        Grid  & $h=1/8$ & $h=1/16$ & $h=1/32$ & $h=1/64$ & $h=1/128$ & $h=1/256$\\
       \hline
       Dim & $d=344$ & $d=1192$ & $d=4424$ & $d=17032$ & $d=66824$ & $d=264712$\\
       \hline \hline
       $\varepsilon_{u_\FF}$ & 7.5836e-4 & 1.6855e-4 & 4.0510e-5 & 1.0011e-5 & 2.4943e-6 & 6.2293e-7\\
       \hline
       $\varepsilon_{v_\FF}$ & 1.5342e-3 & 3.4547e-4 & 8.3952e-5 & 2.0830e-5 & 5.1982e-6 & 1.2991e-6\\
       \hline
       $\varepsilon_{p_\FF}$ & 1.3732e-4 & 3.4712e-5 & 8.6965e-6 & 2.1740e-6 & 5.4331e-7 & 1.3579e-7\\
       \hline 
       $\varepsilon_{p_\PM}$ & 1.9351e-4 & 4.9176e-5 & 1.2384e-5 & 3.1072e-6 & 7.7824e-7 & 1.9474e-7\\
       \hline
    \end{tabular}
    \label{tab:error-2-BJS}
\end{table}
\begin{table}[!ht]
    \caption{Convergence order for all primary variables in Test~\ref{ex:3} (\textit{Problem} $\ABJS$) for physical parameters $\alpha = 1$, $\mu=10^{-3}$ and $k=10^{-2}$}
    \centering
    \begin{tabular}{|c||c|c|c|c|c|}
    \hline
       Grid 1/Grid 2 & $8/16$ & $16/32$ & $32/64$ & $64/128$ & $128/256$ \\
       \hline\hline
    $u_\FF$ &  2.1697 & 2.1133 & 2.0811 & 2.0620 & 2.0499\\
    \hline
    $v_\FF$ & 2.1509 & 2.0959 & 2.0676 & 2.0513 & 2.0412\\
    \hline 
    $p_\FF$ & 1.9840 & 1.9905 & 1.9937 & 1.9954 & 1.9964\\
    \hline
    $p_\PM$ & 1.9764 & 1.9829 & 1.9869 & 1.9895 & 1.9913\\
    \hline
    \end{tabular}
    \label{tab:convergence-rates-ex-2-BJS}
\end{table}

\begin{table}[!ht]
\caption{Errors for all primary variables in Test~\ref{ex:3} (\textit{Problem} $\ABJ$) for physical parameters $\alpha = 1$, $\mu=10^{-3}$ and $k=10^{-2}$
}
    \centering
    \begin{tabular}{|c||c|c|c|c|c|c|}
    \hline
        Grid & $h=1/8$ & $h=1/16$ & $h=1/32$ & $h=1/64$ & $h=1/128$ & $h=1/256$\\
       \hline
       Dim & $d=344$ & $d=1192$ & $d=4424$ & $d=17032$ & $d=66824$ & $d=264712$\\
       \hline \hline
       $\varepsilon_{u_\FF}$ & 9.8945e-4 & 2.1881e-4 & 5.2625e-5 & 1.3012e-5 & 3.2427e-6 & 8.0990e-7\\
       \hline
       $\varepsilon_{v_\FF}$ & 1.6867e-3 & 3.7863e-4 & 9.1928e-5 & 2.2809e-5 & 5.6925e-6 & 1.4227e-6\\
       \hline
       $\varepsilon_{p_\FF}$ & 1.3493e-4 & 3.4003e-5 & 8.5079e-6 & 2.1262e-6 & 5.3137e-7 & 1.3282e-7\\
       \hline 
       $\varepsilon_{p_\PM}$ & 1.9361e-4 & 4.9303e-5 & 1.2428e-5 & 3.1191e-6 & 7.8127e-7 & 1.9550e-7\\
       \hline
    \end{tabular}
    \label{tab:error-2-BJ}
\end{table}

\begin{table}[!ht]
    \caption{Convergence order for all primary variables in Test~\ref{ex:3} (\textit{Problem} $\ABJ$) for physical parameters $\alpha = 1$, $\mu=10^{-3}$ and $k=10^{-2}$}
    \centering
    \begin{tabular}{|c||c|c|c|c|c|}
    \hline
       Grid 1/Grid 2  & $8/16$ & $16/32$ & $32/64$ & $64/128$ & $128/256$ \\
       \hline\hline
    $u_\FF$ & 2.1770 & 2.1164 & 2.0829 & 2.0633 & 2.0509 \\
    \hline
    $v_\FF$ & 2.1553 & 2.0988 & 2.0695& 2.0527 &
    2.0423\\
    \hline 
    $p_\FF$ & 1.9884 & 1.9936 & 1.9959 & 1.9971 &
    1.9977\\
    \hline
    $p_\PM$ & 1.9734 & 1.9807 & 1.9853 & 1.9883 & 1.9903\\
    \hline
    \end{tabular}
    \label{tab:convergence-rates-ex-2-BJ}
\end{table}

\subsection{Eigenvalue distribution}
\begin{figure}[!ht]
    \centering
   \includegraphics[scale = 0.36]{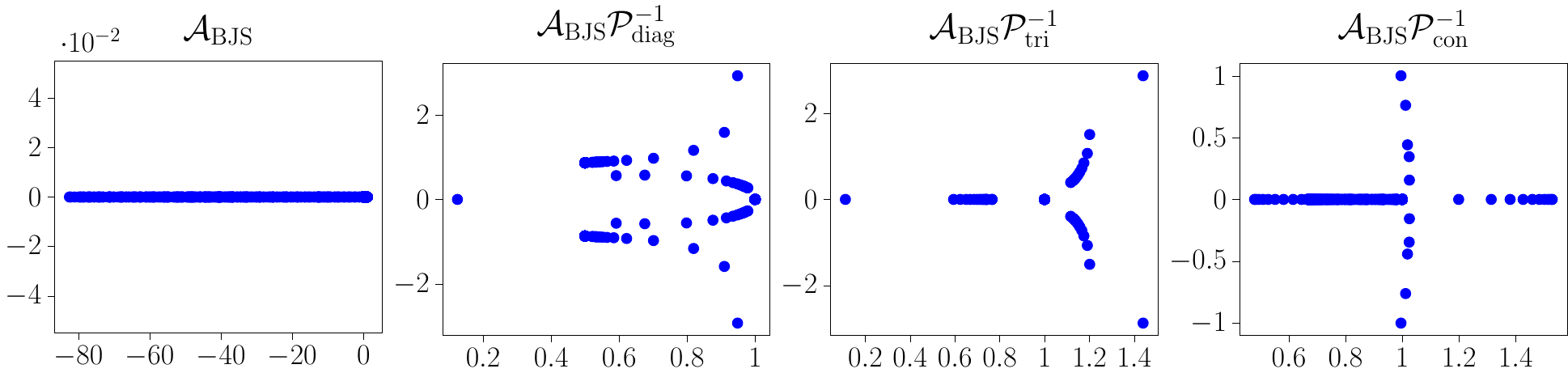}\\[2.5mm]
    \includegraphics[scale = 0.36]{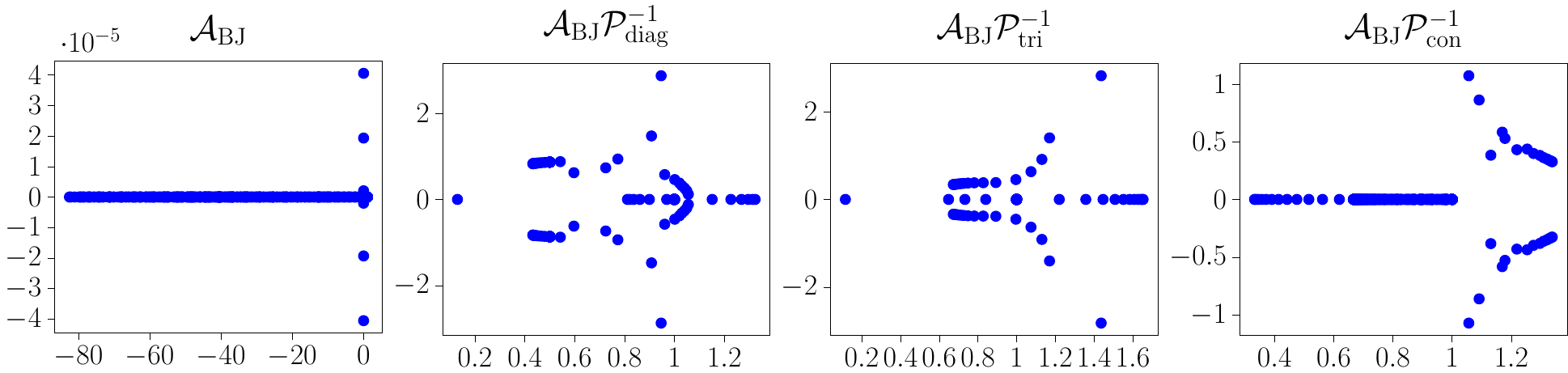}
    \caption{Eigenvalue distributions for the matrices $\mathcal{A}_{\mathrm{BJS}}$ and $\mathcal{A}_{\mathrm{BJ}}$ and the corresponding preconditioned matrices for $h=1/16$ ($d= 1192$)}
    \label{fig:eigenvalues}
\end{figure}

In this section, we first compute and present the eigenvalue distribution for the original and the preconditioned Stokes--Darcy systems. In Fig.~\ref{fig:eigenvalues}, we plot the eigenvalues for the original matrices $\mathcal{A} \in \{\mathcal{A}_{\mathrm{BJ}},\, \mathcal{A}_{\mathrm{BJS}}\}$ given in~\eqref{eq:A-BJS-BJ} and for the corresponding exact preconditioned matrices $\mathcal{A}\mathcal{P}^{-1}$, $\mathcal{P}~\in~\{\Pd,\, \Pt, \, \Pc\}$ from~\eqref{eq:P-D-P-T},~\eqref{eq:P_con}. The physical parameters are chosen as $\alpha =1$, $\mu = 10^{-3}$ and $k = 10^{-2}$. 
As proven in section~\ref{sec:spectral-analysis}, the eigenvalues are clustered and bounded away from zero (Theorems~\ref{theo:eig-APd-2}--\ref{theo:eig-APc}). Therefore, all three developed preconditioners significantly improve the eigenvalue distributions of the original systems~(Fig.~\ref{fig:eigenvalues}).

\begin{figure}[!ht]
    \centering
   \includegraphics[scale = 0.425]{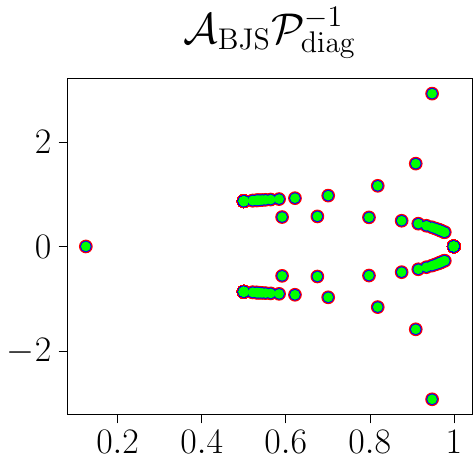}\includegraphics[scale = 0.425]{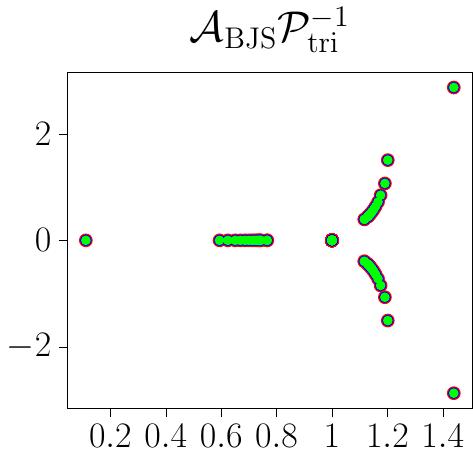}
    \includegraphics[scale = 0.425]{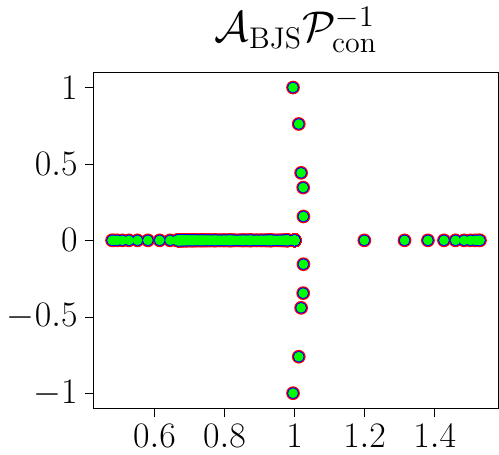}\\[2.5mm]
     \includegraphics[scale = 0.425]{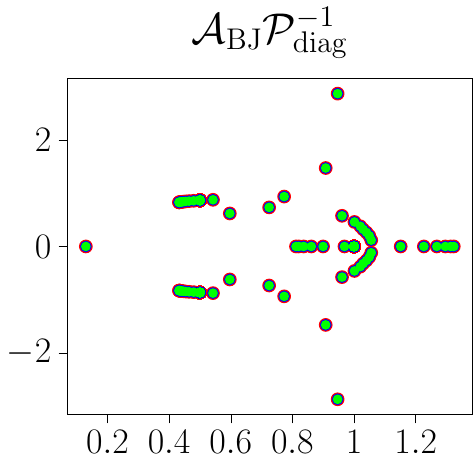}\includegraphics[scale = 0.425]{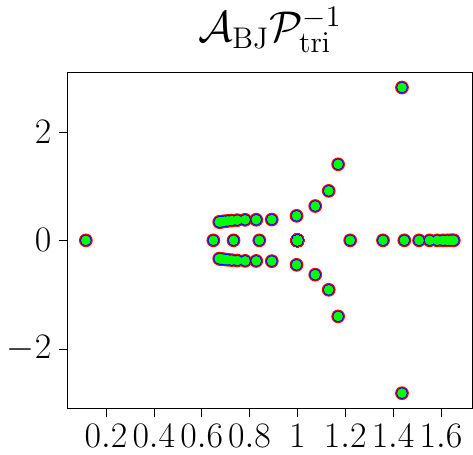}
    \includegraphics[scale = 0.425]{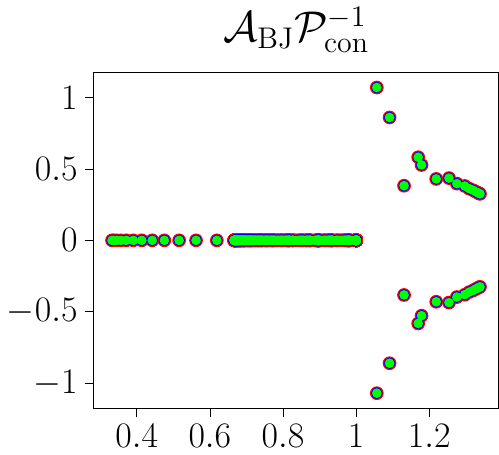}
    \caption{Eigenvalue distributions for the preconditioned matrices for $k = 10^{-2}, \, \alpha = 1$ and $\mu \in \{\mu_1, \mu_2, \mu_3\}$ with $\mu_1 = 10^{-5}$ (\textcolor{red}{$\bullet$}), $\mu_2 = 10^{-3}$ (\textcolor{blue}{$\bullet$}), $\mu_3 = 10^{-1}$ (\textcolor{green}{$\bullet$}) and $h=1/16$ ($d= 1192$)}
    \label{fig:eigenvalues-mu}
\end{figure}

\begin{figure}[!hb]
    \centering
   \includegraphics[scale = 0.425]{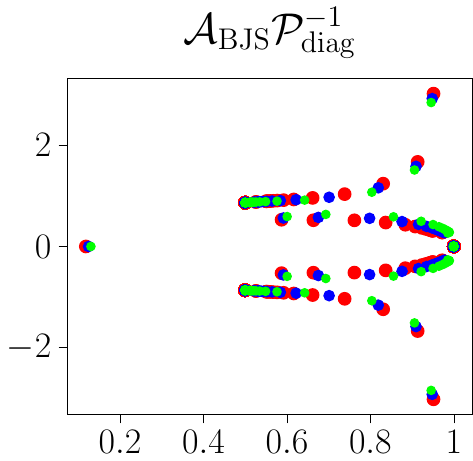}\includegraphics[scale = 0.425]{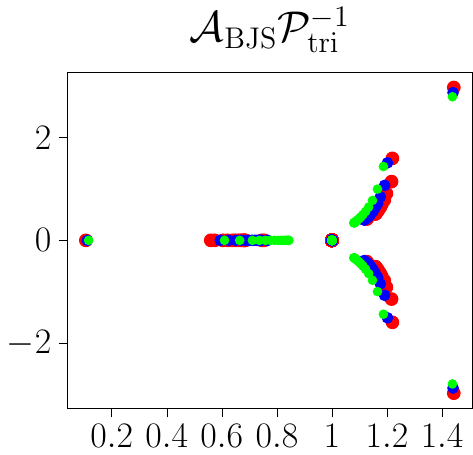}
    \includegraphics[scale = 0.425]{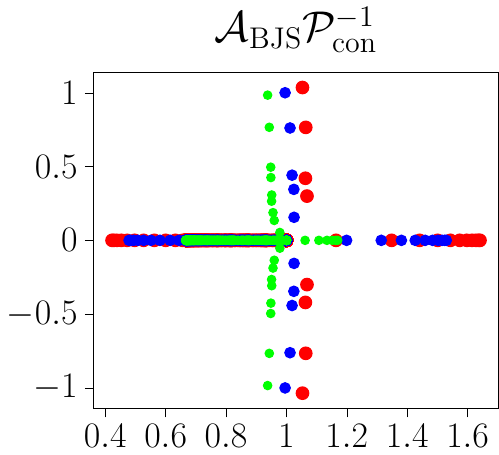}\\[2.5mm]
    \includegraphics[scale = 0.425]{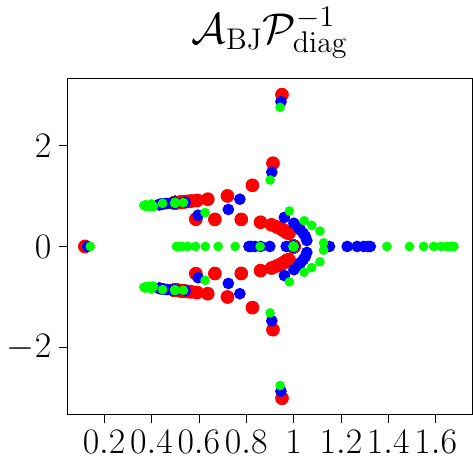}\includegraphics[scale = 0.425]{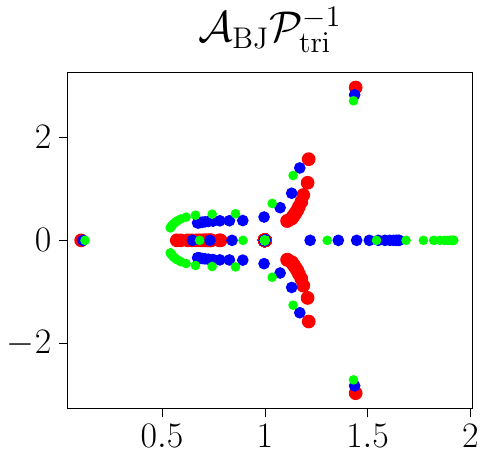}
    \includegraphics[scale = 0.425]{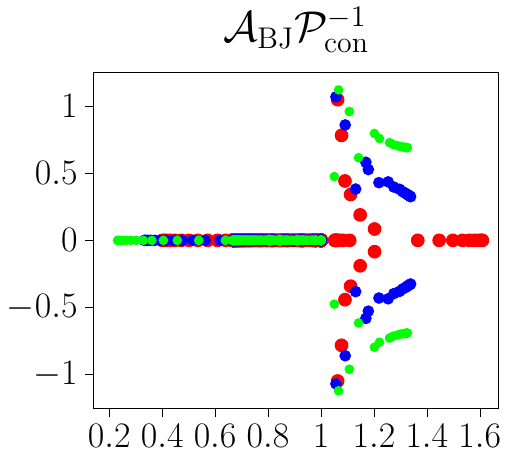}
    \caption{Eigenvalue distributions for the preconditioned matrices for $k = 10^{-2}, \, \mu = 10^{-3}$ and $\alpha \in \{\alpha_1, \alpha_2, \alpha_3\}$ with $\alpha_1 = 10^{-1}$ (\textcolor{red}{$\bullet$}), $\alpha_2 = 1$ (\textcolor{blue}{$\bullet$}), $\alpha_3 = 10$ (\textcolor{green}{$\bullet$}) and $h=1/16$ ($d= 1192$)}
    \label{fig:eigenvalues-alpha}
\end{figure}

\begin{figure}[!hb]
    \centering
   \includegraphics[scale = 0.425]{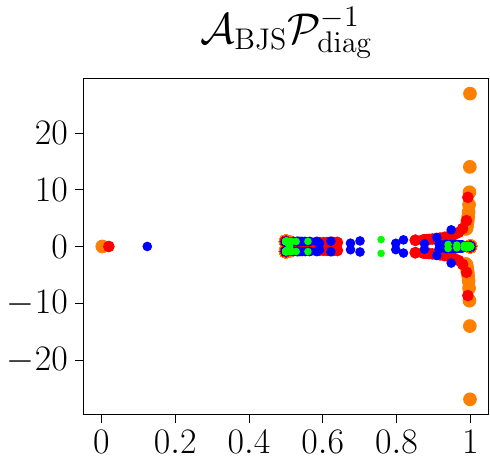}\includegraphics[scale = 0.425]{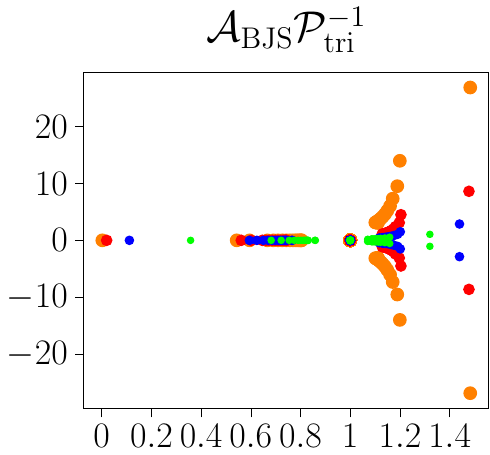}
    \includegraphics[scale = 0.425]{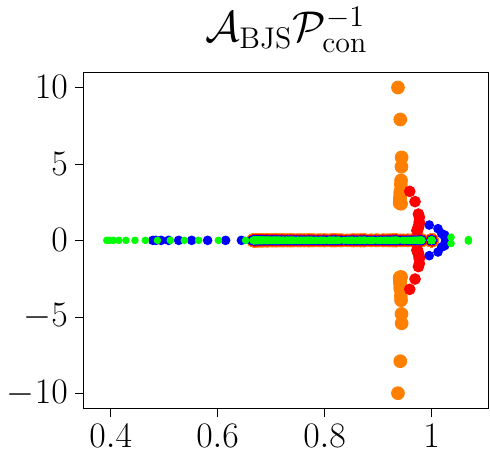}\\[2.5mm]
      \includegraphics[scale = 0.425]{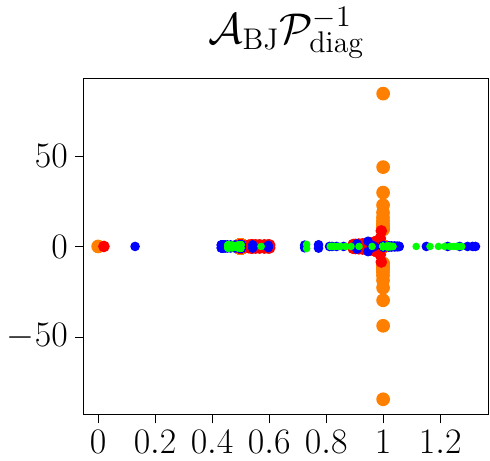}\includegraphics[scale = 0.425]{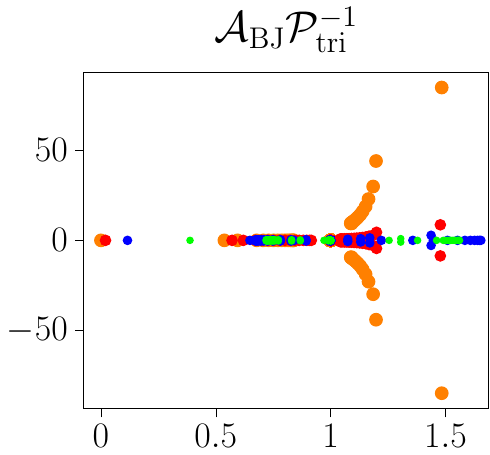}
    \includegraphics[scale = 0.425]{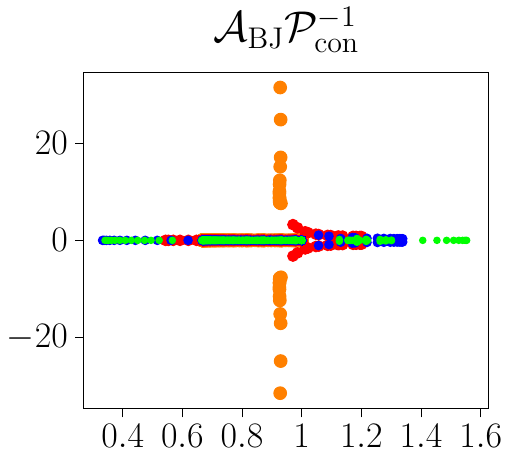}
    \caption{Eigenvalue distributions for the preconditioned matrices for $\mu = 10^{-3}, \, \alpha = 1$ and $k \in \{k_1, k_2, k_3, k_4\}$ with $k_1 = 10^{-5}$ (\textcolor{orange}{$\bullet$}), $k_2 = 10^{-3}$ (\textcolor{red}{$\bullet$}), $k_3 = 10^{-2}$ (\textcolor{blue}{$\bullet$}), $k_4 = 10^{-1}$ (\textcolor{green}{$\bullet$}) and $h=1/16$ ($d= 1192$)}
    \label{fig:eigenvalues-k}
\end{figure}

Moreover, we investigate the influence of the physical parameters on the eigenvalue distributions of the preconditioned matrices. In Fig.~\ref{fig:eigenvalues-mu}, we fix $k=10^{-2}$ and $\alpha=1$, and consider different values for the dynamic viscosity $\mu$. The Beavers--Joseph slip coefficient is varied in Fig.~\ref{fig:eigenvalues-alpha} and various values for the intrinsic permeability are chosen in Fig.~\ref{fig:eigenvalues-k}. 

We observe that the choice of $\mu$ and $\alpha$ have a small impact on the spectra of the preconditioned matrices (Fig.~\ref{fig:eigenvalues-mu}, Fig.~\ref{fig:eigenvalues-alpha}). However, smaller permeability values yield a broader distribution of the eigenvalues (Fig.~\ref{fig:eigenvalues-k}). Nevertheless, for the block diagonal and the block triangular preconditioned matrices, the majority of the eigenvalues (at minimum 95$\%$) are in sphere $B_r$ with radius $r=0.1$ around the clusters identified in Theorems~\ref{theo:eig-APd-2} and~\ref{theo:eig-APt-2} (Tab.~\ref{tab:cluster-eigenvalues-Pd-Pt}). 
Note that, for small permeabilities we obtain eigenvalues close to zero (Fig.~\ref{fig:eigenvalues-k}). However, numerical simulation results show that at most one eigenvalue is in~$B_{0.1}(0)$.

\begin{table}[!ht]
 \caption{Number of eigenvalues in sphere $B_r$ with radius $r=0.1$ around the clusters identified in Theorems~\ref{theo:eig-APd-2} and \ref{theo:eig-APt-2} as well as  around 0 for  $h=1/16$. The total number of eigenvalues is 1192}
    \centering
    \begin{tabular}{|c|c|c||c|c|c||c|c|}
    \hline
    & & &  \multicolumn{3}{c||}{\parbox[0pt][1.5em][c]{0cm}{}$\ABJS\Pd^{-1}$ / $\ABJ\Pd^{-1}$ } & \multicolumn{2}{c|}{\parbox[0pt][1.5em][c]{0cm}{}$\ABJS\Pt^{-1}$ / $\ABJ\Pt^{-1}$}\\
    \hline
    \parbox[0pt][1.5em][c]{0cm}{}$\mu$ & $\alpha$ & $k$ & $B_r(1)$ & $\!\!B_r\!\left((1\pm \sqrt{3}i)/2\right)\!\!$ & $B_r(0)$ & $B_r(1)$ & $B_r(0)$ \\
    \hline
    $10^{-1}$ & $1$ & $10^{-2}$ & 649 / 650 & 502 / 506 &  0 / 0 & 1144 / 1145 &  0 / 0 \\
    $10^{-3}$ & $1$ & $10^{-2}$ & 649 / 650 & 502 / 506 & 0 / 0 & 1144 / 1145 & 0 / 0\\
    $10^{-5}$ & $1$ & $10^{-2}$ & 649 / 650 & 502 / 506 &  0 / 0 & 1144 / 1145 & 0 / 0\\
    \hline
    $10^{-3}$ & $10$ & $10^{-2}$ & 649 / 649 & 504 / 484 & 0 / 0& 1144 / 1144 & 0 / 0\\
    $10^{-3}$ & $1$ & $10^{-2}$ & 649 / 650 & 502 / 506 & 0 / 0 & 1144 / 1145 & 0 / 0\\
    $10^{-3}$ & $10^{-1}$ & $10^{-2}$ & 649 / 649 & 498 / 500 & 0 / 0 & 1144 / 1144 & 0 / 0\\
    \hline
    $10^{-3}$ & $1$ & $10^{-1}$ & 657  / 653 & 510 / 510 & 0 / 0 &  1144 / 1145 & 0 / 0\\
    $10^{-3}$ & $1$ & $10^{-2}$ & 649 / 650 & 502 / 506 & 0 / 0 & 1144 / 1145 & 0 / 0\\
    $10^{-3}$ & $1$ & $10^{-3}$ & 649 / 649 & 480 / 490  & 1 / 1 & 1144 / 1146 & 1 / 1\\
    \hline
    \end{tabular}
    \label{tab:cluster-eigenvalues-Pd-Pt}
\end{table}

In Theorem~\ref{theo:eig-APc}, we identified the clusters $1$ and $\eta= (x^\top A x)/(x^\top G x)$ for the constraint preconditioned matrices. We can bound $\eta$ with the Courant--Fischer theorem as follows
\[\eta_{\min} := \lambda_{\min}\left(\frac{x^\top G^{-1/2}AG^{-1/2}x}{x^\top x}\right) \leq \eta \leq \lambda_{\max}\left(\frac{x^\top G^{-1/2}AG^{-1/2}x}{x^\top x}\right) =: \eta_{\max}.\]
We observe in Tab.~\ref{tab:cluster-eigenvalues-Pc} that almost all eigenvalues ($\geq 98\%$) of the preconditioned matrices are located in the interval $[\eta_{\min}, \eta_{\max}]$. Furthermore, approximately $80\%$ of the eigenvalues are in the sphere $B_{0.1}(1)$. Note that there are no eigenvalues close to zero in the case of the constraint preconditioning.

\begin{table}[!ht]
\caption{Number of eigenvalues around the  clusters identified in Theorem~\ref{theo:eig-APc} for  $h=1/16$. The total number of eigenvalues is $1192$}
    \centering
    \begin{tabular}{|c|c|c||c|c||c|c|}
    \hline
    & & &  & & \multicolumn{2}{c|}{$\ABJS\Pc^{-1}$ / $\ABJ\Pc^{-1}$} \\
    \hline
    \parbox[0pt][1.5em][c]{0cm}{}$\mu$ & $\alpha$ & $k$ & $\eta_{\min}$ & $\eta_{\max}$ & $B_{0.1}(1)$ &  $[\eta_{\min}, \eta_{\max}]$\\
    \hline
    $10^{-1}$ & $1$ & $10^{-2}$ & 3.3325e-1 & 1.6668 & 977 / 962 & 1192 / 1191\\
    $10^{-3}$ & $1$ & $10^{-2}$ & 3.3325e-1 & 1.6668
 &  977 / 962 & 1192 / 1191\\
    $10^{-5}$ & $1$ & $10^{-2}$ & 3.3325e-1 & 1.6668&  977 / 962 & 1192 / 1191\\
    \hline
    $10^{-3}$ & $10$ & $10^{-2}$ & 5.7640e-1
 & 1.4236
 &  980 / 963 & 1192 / 1180\\
    $10^{-3}$ & $1$ & $10^{-2}$ & 3.3325e-1 & 1.6668 &  977 / 962 & 1192 / 1191\\
    $10^{-3}$ & $10^{-1}$ & $10^{-2}$ & 2.4399e-1
 & 1.7560
 &  977 / 973 & 1192 / 1192\\
    \hline
    $10^{-3}$ & $1$ & $10^{-1}$ & 2.6935e-1 & 1.7307 &  978 / 962 & 1192 / 1192\\
    $10^{-3}$ & $1$ & $10^{-2}$ & 3.3325e-1 & 1.6668 &  977 / 962 & 1192 / 1191\\
    $10^{-3}$ & $1$ & $10^{-3}$ & 4.4983e-1 & 1.5502 &  976 / 965 & 1192 / 1192\\
    \hline
    \end{tabular}
    \label{tab:cluster-eigenvalues-Pc}
\end{table}
\subsection{Efficiency analysis}
\begin{figure}[!ht]
    \centering
  \includegraphics[scale = 0.5]{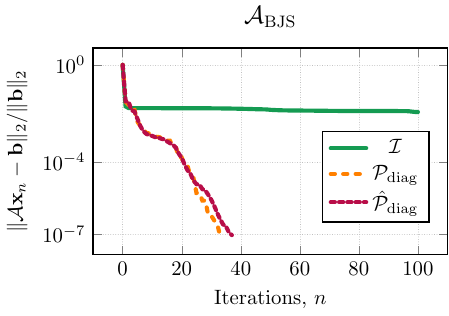}  \hfill\includegraphics[scale = 0.5]{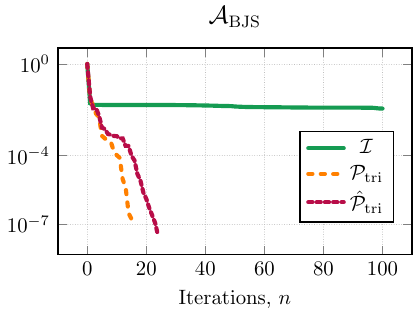}  \hfill
   \includegraphics[scale = 0.5]{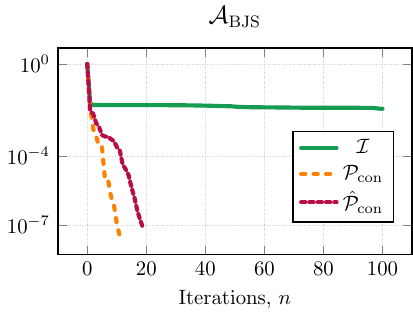} \\
\includegraphics[scale = 0.5]{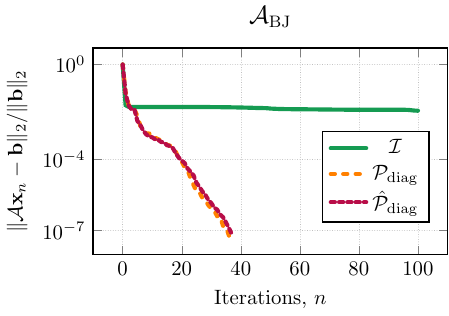}  \hfill\includegraphics[scale = 0.5]{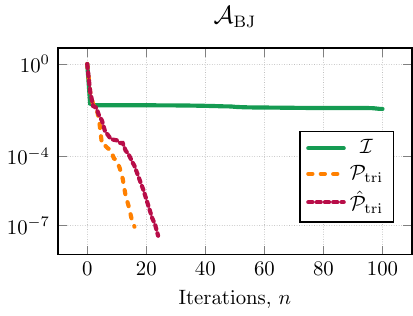}  \hfill
   \includegraphics[scale = 0.5]{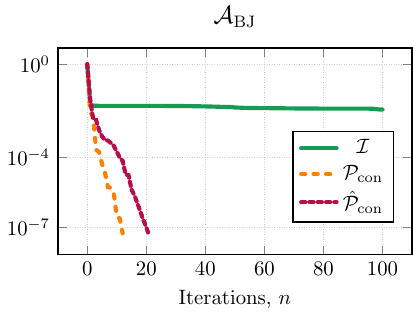}
    \caption{Comparison of different preconditioners for \emph{Problem} $\ABJS$ (top) and \emph{Problem} $\ABJ$ (bottom): original system ($\mathcal{I}$), system with exact preconditioners~($\mathcal{P}$) and system with inexact preconditioners~($\hat{\mathcal{P}}$)  for $h = 1/64$ ($d = 17032$) }
    \label{fig:precon-BJ}
\end{figure}

To study efficiency of the developed preconditioners, we plot the relative residuals \mbox{$\|\mathcal{A}\vec{x}_n-\vec{b}\|_2/\|\vec{b}\|_2$} for $\mathcal{A} \in \{\mathcal{A}_{\mathrm{BJ}},\, \mathcal{A}_{\mathrm{BJS}}\}$ against the number of iterations until the stopping criterion given in the beginning of section~\ref{sec:numerical-results} is reached (Fig.~\ref{fig:precon-BJ}). Here, we consider physical parameters
$\mu = 10^{-3}$, $k = 10^{-2}$ and $\alpha =1$. Note that the GMRES method without preconditioner ($\mathcal{I}$ in Fig.~\ref{fig:precon-BJ}) yields extremely slow convergence. The number of iterations $n$ and the CPU times for the preconditioned GMRES are presented in Tab.~\ref{tab:efficiency-study}. For the inexact versions of the preconditioners $\hat{\mathcal{P}}$, the CPU time is composed of two parts: (i)~time to construct the algebraic grids for $A_{11},\, A_{22}$ and $D$ given in~\eqref{eq:A-BJS-BJ} and~\eqref{eq:P_con}, and (ii)~time to solve the linear system~\eqref{eq:right-preconditioning} using the generated algebraic grids. Note that, even though the exact preconditioners $\Pd$, $\Pt$ and $\Pc$ require less iteration steps (Fig.~\ref{fig:precon-BJ}), the CPU times are significantly higher than for the inexact versions of the preconditioners $\hPd$, $\hPt$ and $\hPc$ (Tab.~\ref{tab:efficiency-study}).

\begin{table}[!ht]
{\footnotesize
\caption{Computational costs to solve \emph{Problems}~$\mathcal{A}_{\mathrm{BJS}}$ and $\mathcal{A}_{\mathrm{BJ}}$ for $h=1/64$ ($d = 17032$)} \label{tab:efficiency-study}
\begin{center}
\begin{tabular}{|c|c|c|c|c|} \hline
& \multicolumn{2}{c|}{\parbox[0pt][1.5em][c]{0cm}{}\textit{Problem} $\mathcal{A}_{\mathrm{BJS}}$} & \multicolumn{2}{c|}{\textit{Problem} $\mathcal{A}_{\mathrm{BJ}}$} \\
\hline
\parbox[0pt][1.5em][c]{0cm}{}Preconditioner & Iterations & CPU time [s] & Iterations & CPU time [s] \\
     \hline
     $\Pd$ & 33 & 31.14 & 36 & 35.66\\
     $\hPd$ & 37 & 0.5103 + 0.2311 & 37 & 0.6347 + 0.2614\\
     \hline
     $\Pt$ & 16 & 40.86  & 16 & 40.83 \\
     $\hPt$ & 23 & 0.4769 + 0.1598 & 24& 0.4606 + 0.1609 \\
     \hline
     $\Pc$ & 11 & 34.23 & 12 & 36.61\\
     $\hPc$ & 18& 0.4934 + 0.1263 & 21 &0.5348 +  0.1530 \\
     \hline
\end{tabular}
\end{center}
}
\end{table}

\subsection{Robustness analysis}
An important property of preconditioners is the independence of the convergence rate of the iterative method from the grid width $h$. As above, we fix the parameters $\mu = 10^{-3}$, $k = 10^{-2}$, $\alpha = 1$ and study the convergence for different grid widths~$h$. Due to the large size of the linear systems, we provide the results only for the inexact versions of the preconditioners. As it can be seen in Tab.~\ref{tab:h-analysis}, the number of iteration steps stays nearly constant. For $h\to 0$ the number of iterations slightly decreases. This is due to the fact that the contributions of the blocks $C$ and $C_1$ are of order $O(h)$, yielding $\ABJS \rightarrow \Abar$ for~$h\to 0$. We observe the same behaviour for $\ABJ$ under the additional condition $\alpha k_{yy}/\sqrt{k_{xx}} \ll 1$, which also appears in Theorems~\ref{theo:eig-APd-2} and~\ref{theo:eig-APt-2}. 
\begin{table}[tbhp]
\caption{Number of iteration steps to solve \emph{Problems}~$\mathcal{A}_{\mathrm{BJS}}$ and $\mathcal{A}_{\mathrm{BJ}}$ for $\mu = 10^{-3}$, $\alpha = 1$ and $k = 10^{-2}$ and different grid widths~$h$ (dimensions of linear systems~$d$)}
    \label{tab:h-analysis}
\begin{center}
\begin{tabular}{|C{0.9cm}|C{1.1cm}||C{0.9cm}|C{0.9cm}|C{0.9cm}|C{0.9cm}|C{0.9cm}|C{0.9cm}|}
\hline
        & & \multicolumn{3}{c|}{\parbox[0pt][1.5em][c]{0cm}{}\emph{Problem}~$\mathcal{A}_{\mathrm{BJS}}$} & \multicolumn{3}{c|}{\emph{Problem}~$\mathcal{A}_{\mathrm{BJ}}$} \\
       \hline
       \parbox[0pt][1.5em][c]{0cm}{} $h$ & $d$ & $\hPd$ & $\hPt$ & $\hPc$ & $\hPd$ & $\hPt$ & $\hPc$ \\
        \hline
   \parbox[0pt][1.5em][c]{0cm}{} $2^{-3}$ & 344 & 37 & 26 & 21 & 39 & 27 & 23 \\
\hline
\parbox[0pt][1.5em][c]{0cm}{} $2^{-4}$ & 1192 & 39 & 25 & 20 & 41 & 27 & 24\\
\hline
\parbox[0pt][1.5em][c]{0cm}{} $2^{-5}$ & 4424 & 38 & 24 & 20 &40 & 25 & 22\\
\hline
\parbox[0pt][1.5em][c]{0cm}{} $2^{-6}$ & 17032 & 37 & 23 & 19 & 37 & 24 & 21\\
\hline
\parbox[0pt][1.5em][c]{0cm}{} $2^{-7}$ & 66824 & 35 & 22 & 18 & 35 & 22 & 19\\
\hline
\parbox[0pt][1.5em][c]{0cm}{} $2^{-8}$ & 264712 & 32 & 21 & 17 & 32 & 21 & 17\\
\hline
\parbox[0pt][1.5em][c]{0cm}{} $2^{-9}$ & 1053704 & 29 & 20 & 16 & 29 & 20 & 17\\
\hline
    \end{tabular}
\end{center}
\end{table} 

\begin{table}[!ht]
\caption{Number of iteration steps to solve \emph{Problems}~$\mathcal{A}_{\mathrm{BJS}}$ and $\mathcal{A}_{\mathrm{BJ}}$ for different values of parameters $\mu$, $\alpha$ and $k$ for~$h=1/64$ ($d = 17032$) }
    \centering
    \begin{tabular}{|c|c|c||c|c|c|c|c|c|}
    \hline
    & & & \multicolumn{3}{c|}{\parbox[0pt][1.5em][c]{0cm}{} \textit{Problem} $\mathcal{A}_{\mathrm{BJS}}$} & \multicolumn{3}{c|}{\textit{Problem} $\mathcal{A}_{\mathrm{BJ}}$} \\
    \hline
    $\mu$ & $\alpha$ & $k$ & $\hPd$ & $\hPt$ & $\hPc$ & $\hPd$ & $\hPt$ & $\hPc$\\
    \hline
    $10^{-1}$ & $1$ & $10^{-2}$ & 42 & 26 & 22 & 44 & 26 & 23\\
    $10^{-2}$ & $1$  & $10^{-2}$ & 38 & 23 & 20 & 39 & 24 & 21\\
    $10^{-3}$ & $1$ & $10^{-2}$ & 37 & 23 & 19 & 37 & 24 & 21\\
    $10^{-4}$ & $1$ & $10^{-2}$ & 37 & 23 &19 & 37 & 24 & 21\\
    $10^{-5}$ & $1$ & $10^{-2}$ & 37 & 23 & 19 & 37 & 24 & 21\\
    \hline
    $10^{-3}$ & $10$ & $10^{-2}$ & 33 & 24 & 19 & 40 & 27 & 25\\
    $10^{-3}$ & $1$ & $10^{-2}$ & 37 & 23 & 19 & 37 & 24 & 21\\
    $10^{-3}$ & $10^{-1}$ & $10^{-2}$ & 39 & 23 & 20 & 39 & 23 & 19\\
    \hline
    $10^{-3}$ & $1$ & $10^{-2}$ & 37 & 23 & 19 & 37 & 24 & 21\\
    $10^{-3}$ & $1$ & $10^{-3}$ & 53 & 38 & 32 & 54 & 38 & 35\\
    $10^{-3}$ & $1$ & $10^{-4}$ & 84 & 67 & 60 & 84 & 67 & 60\\
    $10^{-3}$ & $1$ & $10^{-5}$ & 146 & 121 & 105 & 145 & 120 & 105\\
    $10^{-3}$ & $1$ & $10^{-8}$ & 155 & 140 & 116 & 155 & 140 & 116\\
    \hline
    \end{tabular}
    \label{tab:robustness-analysis}
\end{table}

To study robustness of the preconditioners with respect to the physical parameters $\mu,\,k$ and $\alpha$, we consider different values of the intrinsic permeability $k$, the dynamic viscosity $\mu$ and the Beavers--Joseph slip parameter $\alpha$. 
We observe high robustness with respect to viscosity $\mu$ and Beavers--Joseph parameter $\alpha$ (Tab.~\ref{tab:robustness-analysis}) which is in correspondence with the results on the eigenvalue distributions (Fig.~\ref{fig:eigenvalues-mu}, Fig.~\ref{fig:eigenvalues-alpha}).
The robustness of the preconditioners regarding the permeability~$k$ is moderate. For smaller permeability values the number of iteration steps is increasing (Tab.~\ref{tab:robustness-analysis}). The same behaviour was observed in the literature for other types of preconditioners for coupled Stokes--Darcy problems, e.g. \cite{Greif_He_2023, Holter-et-al-21, Boon-et-al-22}.

\section{Conclusions}
\label{sec:conclusions}
In this paper, we proposed and analysed three different preconditioners for coupled Stokes--Darcy systems: a block diagonal, a block triangular and a constraint preconditioner. We considered two sets of interface conditions with either the Beavers--Joseph--Saffman (\emph{Problem} $\ABJS$) or the original Beavers--Joseph coupling condition (\emph{Problem} $\ABJ$). We developed the second order finite volume method on staggered grids (MAC scheme) to discretise the coupled Stokes--Darcy problems paying special attention on the discretisation of coupling conditions. We used the right-preconditioned GMRES method to solve the resulting linear systems. 

We analysed the spectra identifying the clusters for the eigenvalues for the Beavers--Joseph--Saffman and the more general Beavers--Joseph case. In addition, for the Beavers--Joseph--Saffman condition, we provided bounds on the spectrum for the exact variants of the developed preconditioners.
To confirm the obtained theoretical results, we performed a series of numerical experiments for Stokes--Darcy problems. The numerical simulation results show that both, the exact and inexact variants of the preconditioners, significantly improve convergence of the GMRES method. Even though the number of iteration steps is smaller for the exact variants of the preconditioners, the inexact versions provide remarkably smaller CPU times. 

The preconditioners are robust with respect to varying grid width $h$. 
We also observe robustness for different values of the viscosity $\mu$ and the Beavers--Joseph slip coefficient $\alpha$.
However, with decreasing intrinsic permeability $k$ the number of iteration steps increases. Such behaviour was already observed in the literature for different types of preconditioners for Stokes--Darcy problems. 

To guarantee convergence of the GMRES method with the proposed preconditioners, the field-of-values analysis is needed. This analysis is subject of the consequent manuscript.

\section*{Funding}
The work is funded by the Deutsche Forschungsgemeinschaft (DFG, German Research Foundation) – Project Number 327154368 – SFB 1313 and Project Number 490872182.

\section*{Data availability}
The datasets generated and analysed during the current study are available from the corresponding author on reasonable request.

\section*{Declarations}
\textbf{Conflict of interest} The authors declare that they have no conflict of interest.

\bibliographystyle{spmpsci}      
\bibliography{references} 

\end{document}